\documentclass[11pt]{article}
\usepackage{amssymb}
\usepackage{amsmath}
\usepackage{mathrsfs}
\usepackage{graphics}
\usepackage{graphicx}
\usepackage{xcolor}
\usepackage{subfigure}
\usepackage[T1]{fontenc}
\usepackage{latexsym,amssymb,amsmath,amsfonts,amsthm}\usepackage{txfonts}
\newcommand{\be}{\begin{eqnarray}}
\newcommand{\ben}{\begin{eqnarray*}}
\newcommand{\en}{\end{eqnarray}}
\newcommand{\enn}{\end{eqnarray*}}

\newtheorem{theorem}{Theorem}[section]
\newtheorem{lemma}{Lemma}[section]
\newtheorem{prp}[theorem]{Proposition}
\newtheorem{thm}[theorem]{Theorem}

\newtheorem{dfn}{Definition}[section]

\begin{document}
\renewcommand{\theequation}{\arabic{section}.\arabic{equation}}
\begin{titlepage}
\title{\bf A moderate deviation principle for 2D
stochastic primitive equations\thanks{This work was partially
supported by NNSF of China(Grant No. 11401057),   Natural Science Foundation Project of CQ  (Grant No. cstc2016jcyjA0326),
Fundamental Research Funds for the Central Universities(Grant No. 106112015CDJXY100005) and China Scholarship Council (Grant No.201506055003).}}

\author{Rangrang Zhang$^{1,}$ Guoli Zhou $^{2,}$\\
{\small $^1$ Department of  Mathematics,
Beijing Institute of Technology,}\\
{\small No. 5 Zhongguancun South
Street, Haidian District, Beijing, 100081, P R China}\\
{\small $^2$ School of Statistics and Mathematics, Chongqing University,}\\
 {\small No.174 Shazheng Street, Shapingba District, Chongqing, 400044, P R China}\\
({\sf rrzhang@amss.ac.cn}\ {\sf zhouguoli736@126.com})}
\date{}
\end{titlepage}
\maketitle

\textbf{Abstract}:
 In this paper, we establish a central limit theorem and a moderate deviations for 2D stochastic primitive equations with multiplicative noise. The proof is mainly based on the weak convergence approach.



\section{Introduction}
As a fundamental model in meteorology, the primitive equations
were derived from the Navier-Stokes equations, with rotation, coupled with thermodynamics
and salinity diffusion-transport equations, by assuming two important simplifications: the Boussinesq approximation and the hydrostatic balance (see \cite{L-T-W-1,L-T-W-2,JP} and the references therein). This model in the deterministic case has been intensively investigated because of the interests stemmed from physics and mathematics. For example, the mathematical study of the primitive equations originated in a series of articles by Lions, Temam, and Wang in the early 1990s (see \cite{L-T-W-1,L-T-W-2,L-T-W-3,L-T-W-4}), where they set up the mathematical framework and showed the global existence of weak solutions. Cao and Titi developed a beautiful approach to dealing with the $L^6$-norm of the fluctuation $\tilde{v}$ of horizontal velocity and obtained the global well-posedness for the 3D viscous primitive equations in  \cite{C-T-1}.
\

For the primitive equations in random case, many results have been obtained. Debussche, Glatt-Holtz, Temam and Ziane established the global well-posedness of the strong solution of the primitive
equations driven by multiplicative random noises in \cite{D-G-T-Z}. The ergodic theory of 3D stochastic primitive equations driven by regular multiplicative noise was studied in \cite{RR}, where we proved that all weak solutions which are limits of spectral Galerkin approximations share the same invariant measure. Using a new method we proved the existence of random attractor for 3D stochastic primitive equations driven by fractional noise in $\cite{Zgl}.$
In \cite{G-S}, Gao and Sun obtained a Freidlin-Wentzell's large deviation principle (LDP) for the stochastic primitive equations in two dimensional case while we established the same result in three dimensional case in \cite{R-R}.

In this paper, we shall investigate deviations of the strong solution $Y^\varepsilon$ (see (\ref{e-1}))of 2D stochastic primitive equations from the deterministic solution $Y^0$ (see (\ref{e-2})), as $\varepsilon$ decreases to 0, that is, the asymptotic behavior of the trajectory,
\[
Z^\varepsilon(t)=\frac{1}{\sqrt{\varepsilon}\lambda(\varepsilon)}(Y^\varepsilon-Y^0)(t), \quad t\in[0,T],
\]
where $\lambda(\varepsilon)$ is some deviation scale which strongly influences the asymptotic behavior of $Z^\varepsilon$. Concretely, three cases are involved.
\begin{description}
  \item[(1)] The case $\lambda(\varepsilon)=\frac{1}{\sqrt{\varepsilon}}$ provides some large deviation principles (LDP), which has been obtained by Gao et al. in \cite{G-S}.
  \item[(2)] The case $\lambda(\varepsilon)=1$  provides the central limit theorem (CLT). We will show that $Z^\varepsilon$ converges to a solution of a stochastic equation as $\varepsilon$ decrease to 0 in Sect. \ref{CLT}.
  \item[(3)] To fill in the gap between the CLT scale ($\lambda(\varepsilon)=1$) and the large deviations scale ($\lambda(\varepsilon)=\frac{1}{\sqrt{\varepsilon}}$), we will study the so-called moderate deviation principle (MDP) in Sect. \ref{MDP}. Here, the deviations scale satisfies
      \begin{eqnarray}\label{e-43}
      \lambda(\varepsilon)\rightarrow +\infty,\ \sqrt{\varepsilon}\lambda(\varepsilon)\rightarrow 0\quad {\rm{as}} \ \varepsilon\rightarrow 0.
      \end{eqnarray}
\end{description}

Similar to LDP, MDP arises in the theory of statistical inference naturally, which can provide us with the rate of convergence and a useful method for constructing asymptotic confidence intervals ( see, e.g. \cite{E-1,E-2,I-K,K} and references therein).
The proof of moderate deviations is mainly based on the weak convergence approach, which is developed by Dupuis and Ellis in \cite{DE}. The key idea is to prove some variational representation formula about
the Laplace transform of bounded continuous functionals, which will lead to proving an equivalence between
the Laplace principle and LDP. In particular, for Brownian functionals, an elegant variational
representation formula has been established by Bou\'{e}, Dupuis \cite{MP} and Budhiraja, Dupuis \cite{BD}.

Up to now, there are some results about the moderate deviations for fluid dynamics models and other processes. For example, Wang, Zhai and Zhang \cite{WZZ} established the CLT and MDP for 2D Navier-Stokes equations with multiplicative Gaussian noise in the state space $C([0,T];H)\cap L^2([0,T];V)$. Further, the MDP for 2D Navier-Stokes equations driven by multiplicative L\'{e}vy noises were considered by Dong et al. \cite{DXZZ} in state space $D([0,T];H)\cap L^2([0,T];V)$.
In view of the characterization of the super-Brownian motion (SBM) and the Fleming-Viot process (FVP), Fatheddin and Xiong obtained the MDP for those processes in \cite{FX}.
In this paper, we consider the CLT and MDP of 2D stochastic primitive equations, the state space is chosen to be $C([0,T];H)\cap L^2([0,T]; V)$.

Compared with 2D Navier-Stokes equations, the 2D primitive equations are more complex because the nonlinear term of the 2D NSEs is of the form velocity
$\times$ derivative of velocity, while the 2D viscous PEs has a challenging term:
 derivative of horizontal velocity $\times$ derivative of horizontal velocity which results in essential difficulties in obtaining moment estimates. An additional term $|\partial_z Y|^{2}$ appears when estimating the nonlinear term $ B(Y,Y) $ in the functional space $H$. To overcome this difficulty, we take advantage of the special geometrical structure of 2D primitive equations to obtain estimate of $|\partial_z Y|^{2}$, which is also essential to establish some tightness results. With the help of those estimates, a central limits theorem and a moderate derivations principle for 2D primitive equations are proved.


This paper is organized as follows. The mathematical framework of 2D primitive equations is in Sect. 2. In Sect. 3, the formulation of 2D primitive equations: functional spaces, hypothesises on the noise and definition of solution are presented. The central limit theorem is proved in Sect. 4. We obtain the moderate deviation principle by using the weak convergence method in Sect. 5.
%

\section{Preliminaries}
Consider the
2D primitive equations driven by a stochastic forcing in a Cartesian system
\begin{eqnarray}\label{eq-1}
\frac{\partial v}{\partial t}-\mu_1\Delta v+v\partial_{x} v+\theta\partial_{z} v+\partial_{x} P &=&\psi_1(t,v,T)\frac{dW_1}{dt},\\
\label{eq-2}
\partial_{z}P+T&=&0,\\
\label{eq-3}
\partial_{x} v+\partial_{z}\theta&=&0, \\
\label{eq-4}
\frac{\partial T}{\partial t}-\mu_2\Delta T+v\partial_{x} T+\theta\partial_{z} T&=&\psi_2(t,v,T)\frac{dW_2}{dt},
\end{eqnarray}
where the velocity $v=v(t,x,z)\in \mathbb{R}$, the vertical velocity $\theta$, the temperature\ $T$ and the pressure\ $P$ are all unknown functionals. $(x,z)\in \mathcal{M}=(0,L)\times (-h,0)$, $t>0$, $W_1$ and $W_2$ are two independent cylindrical Wiener processes, which will be given in Sect. \ref{sec 2}. $\Delta=\partial^{2}_{x}+\partial^{2}_{z}$ is the Laplacian operator.  Without loss of generality, we assume that
$$
\mu_1=\mu_2=1.
$$
We impose the following boundary conditions:
\begin{eqnarray}\label{eq1}
 \partial_{z}v=0, \ \theta=0,\  \partial_{z}T=0  & &  \rm {on} \ \Gamma_{\textit{u}}=(0,L)\times \{0\},\\
 \label{eq2}
 \partial_{z}v=0,\  \theta=0,\  \partial_{z}T=0        & &   \rm {on}\ \Gamma_{\textit{b}}=(0,L)\times\{-h\},\\
 \label{eq3}
 v=0,\ \partial_{x} T=0  & &   \rm {on}\  \Gamma_{l}=\{0,L\}\times (-h,0).
 \end{eqnarray}
Integrating (\ref{eq-3}) from $-h$ to $z$ and using (\ref{eq1})-(\ref{eq2}), we have
\begin{equation}
\theta(t,x,z):=\Phi(v)(t,x,z)=-\int^{z}_{-h}\partial_x v (t,x,z')dz',
\end{equation}
moreover, in view of (\ref{eq-3}) and (\ref{eq1})-(\ref{eq3}),
\[
\int^{0}_{-h} v dz=0.
\]
Integrating (\ref{eq-2}) from $-h$ to $z$, set $p_{b}$ be a certain unknown function at $\Gamma_{b}$ satisfying
\[
P(x,z,t)= p_{b}(x,t)-\int^{z}_{-h} T(x,z',t) dz'.
\]
Then, (\ref{eq-1})-(\ref{eq-4}) can be rewritten as
 \begin{eqnarray}\label{eq5-1}
&\frac{\partial v}{\partial t}-\Delta v+v\partial_{x} v+\Phi(v)\partial_{z} v+\partial_{x} p_b-\int^z_{-h}\partial_x T dz' =\psi_1(t,v,T)\frac{dW_1}{dt},&\\
\label{eq-6-1}
&\frac{\partial T}{\partial t}-\Delta T+v\partial_{x} T+\Phi(v)\partial_{z} T=\psi_2(t,v,T)\frac{dW_2}{dt},&\\
\label{eq-7-1}
&\int^{0}_{-h} v dz=0.&
\end{eqnarray}
The boundary value conditions for (\ref{eq5-1})-(\ref{eq-7-1}) are given by
\begin{eqnarray}\label{eq-8-1}
 \partial_{z}v=0,\ \partial_{z}T=0    &&      {\rm on}\ \Gamma_{\textit{u}},\\
 \label{eq-9-1}
 \partial_{z}v=0,\ \partial_{z}T=0     &&   {\rm on}\ \Gamma_{\textit{b}},\\
 \label{eq-10-1}
 v=0,\ \partial_{x} T=0      &&   {\rm on}\ \Gamma_{l}.
 \end{eqnarray}
Denote $Y =(v, T)$ and the initial condition
 \begin{equation}\label{eq-11-1}
 Y(0)=Y_0=(v_{{0}},T_{{0}}).
 \end{equation}
\section{Formulation of the SPDE}\label{sec 2}
\subsection{Some Functional Spaces }
Let $\mathcal{L}(K_1;K_2)$ (resp. $\mathcal{L}_2(K_1;K_2)$) be the space of bounded (resp. Hilbert-Schmidt) linear operators from the Hilbert space $K_1$ to $K_2$, whose norm is denoted by $\|\cdot\|_{\mathcal{L}(K_1;K_2)}(\|\cdot\|_{\mathcal{L}_2(K_1;K_2)})$. For $p\in \mathbb{Z}^+$, set
\begin{eqnarray*}
 |\phi|_p=\left\{
            \begin{array}{ll}
              \left(\int_{\mathcal{M}}|\phi(x,z)|^pdxdz\right)^{\frac{1}{p}}, & \phi\in L^p(\mathcal{M}), \\
               \left(\int^L_{0}|\phi(x)|^pdx\right)^{\frac{1}{p}}, &  \phi\in L^p((0,L)).
            \end{array}
          \right.
\end{eqnarray*}
In particular, $|\cdot|$ and $(\cdot,\cdot)$ represent norm and inner product of $L^2(\mathcal{M})$ (or $L^2((0,L))$), respectively. For $m\in \mathbb{N}_+$, $(W^{m,p}(\mathcal{M}), \|\cdot\|_{m,p})$ stands for the classical Sobolev space, see \cite{Adams}. When $p=2$, we denote by $H^m(\mathcal{M})=W^{m,2}(\mathcal{M})$,
\begin{equation}\notag
\left\{
  \begin{array}{ll}
    H^{m}(\mathcal{M})=\Big\{Y\Big| \partial_{\alpha}Y\in (L^2(\mathcal{M}))^2\ {\rm for} \ |\alpha|\leq m\Big\},&  \\
    |Y|^2_{H^{m}(\mathcal{M})}=\sum_{0\leq|\alpha|\leq m}|\partial_{\alpha}Y|^2. &
  \end{array}
\right.
\end{equation}
It's known that $(H^{m}(\mathcal{M}), |\cdot|_{H^{m}(\mathcal{M})})$ is a Hilbert space. $|\cdot|_{H^p((0,L))}$ stands for the norm of $H^p((0,L))$ for $p\in \mathbb{Z}^+$.

 Define working spaces for (\ref{eq5-1})-(\ref{eq-11-1}):
 \begin{eqnarray}\notag
 &&\bar{H}:=\left\{(v,T)\in (C^{\infty}(\mathcal{\bar{M}}))^{2}:\  (v,T)\ \mathrm{satisfies}\  \mathrm{boundary}\ \mathrm{conditions}\  (\ref{eq-8-1})-(\ref{eq-10-1}) \right\}, \notag
 \end{eqnarray}
Let $H$ and $V$ to be the closure of $ \bar{H}$ under the topology of $L^{2}(\mathcal{M})$ and $H^{1}(\mathcal{M})$ respectively. Inner product in $H$ is
\begin{eqnarray*}
(Y,\tilde{Y})=(v,\tilde{v})+(T,\tilde{T})=\int_{\mathcal{M}}(v\tilde{v}+T\tilde{T})dxdz,\quad
|Y|=(Y,Y)^{\frac{1}{2}}.
\end{eqnarray*}
Taking into account the boundary conditions (\ref{eq-8-1})-(\ref{eq-10-1}), the inner product and norm on $V$ can be given by
\begin{eqnarray*}
((Y,\tilde{Y}))&=&((v,\tilde{v}))_1+((T,\tilde{T}))_2,\\
((v,\tilde{v}))_1&=&\int_{\mathcal{M}}(\partial_x v  \partial_x \tilde{v}+\partial_z v  \partial_z \tilde{v})dxdz,\\
((T,\tilde{T}))_2&=&\int_{\mathcal{M}}(\partial_x T \partial_x \tilde{T}+\partial_z T  \partial_z \tilde{T})dxdz,
\end{eqnarray*}
and take $\|\cdot\|=\sqrt{((\cdot,\cdot))}$, where $Y=(v,T), \tilde{Y}=(\tilde{v},\tilde{T})\in V$.
Note that under the above definition, a Poincar$\acute{e}$ inequality $|Y|\leq C\|Y\|$ holds for all $Y\in V$. Let $\check{V}$ be the closure of $V\cap (C^{\infty}(\bar{\mathcal{M}}))^2$ in $(H^{2}(\mathcal{M}))^2$
and equip this space with the norm and inner product of $H^{2}(\mathcal{M})$.

Define the intermediate space
\[
\mathcal{H}=\{Y=(v,T)\in H,\ \partial_z Y=(\partial_z v,\partial_z T )\in H\}.
\]

Let $V'$ be the dual space of $V$. Then, the dense and continuous embeddings
\[
V\hookrightarrow H = H'\hookrightarrow V',
\]
hold and denote by $\langle x, y\rangle$ the duality between $x\in V$ and $y\in V'$.
\subsection{Some Functionals}
Define $P_H$ be the Leray type projection operator from $(L^{2}(\mathcal{M}))^2$ onto $H$.
The principle linear portion of the equation is defined by
\begin{eqnarray*}
AY:=P_H \left(                 
  \begin{array}{c}   
   -\Delta v\\  
   -\Delta T \\  
  \end{array}
\right) \quad {\rm{for}} \ Y=(v,T)\in D(A)
\end{eqnarray*}
where
\begin{eqnarray*}
D(A)\ \mathrm{is}\ \mathrm{the}\ \mathrm{closure}\ \mathrm{of}\ \bar{H}\ \mathrm{with}\ \mathrm{respect}\ \mathrm{to}\ \mathrm{the}\ \mathrm{topology}\ \mathrm{of}\ H^{2}(\mathcal{M}).
\end{eqnarray*}
It's well-known that $A$ is a self-adjoint and positive definite operator. Due to the regularity results of the Stokes problem of geophysical fluid dynamics, we have $|AY|\cong |Y|_{H^2(\mathcal{O})}$ ( see \cite{Z}).

For $Y=(v,T)$, $\tilde{Y}=(\tilde{v},\tilde{T})\in D(A)$, define $B(Y, \tilde{Y})=B_1(v, \tilde{Y})+B_2(v, \tilde{Y})$, where
\begin{eqnarray*}
B_1(v, \tilde{Y}):=P_H \left(                 
  \begin{array}{c}   
   v\partial_x\tilde{v}\\  
    v\partial_x \tilde{T} \\  
  \end{array}
\right), \quad
B_2(v, \tilde{Y}):=P_H \left(                 
  \begin{array}{c}   
   \Phi(v)\partial_z \tilde{v}\\  
   \Phi(v)\partial_z \tilde{T} \\  
  \end{array}
\right)
\end{eqnarray*}
By interpolation inequalities (see \cite{G-S, GHZ08}, etc), we have
\begin{lemma}\label{lemma-1}
 For any $Y=(v,T),\ \tilde{Y}, \hat{Y}\in V$, there exists a constant $C$ such that
 \begin{eqnarray}\label{ee-1}
  \langle v\partial_x{v}+\Phi(v)\partial_z {v}, \partial_{zz} v\rangle=0.\\
  \label{ee-2}
 \langle B(Y,\tilde{Y}),\hat{Y}\rangle=-\langle B(Y,\hat{Y}),\tilde{Y}\rangle, \quad \langle B(Y,\tilde{Y}),\tilde{Y}\rangle=0.\\
 \label{ee-3}
|\langle B(Y,\tilde{Y}), \hat{Y}\rangle|\leq C\|\tilde{Y}\||Y|^{\frac{1}{2}}\|Y\|^{\frac{1}{2}}|\hat{Y}|^{\frac{1}{2}}\|\hat{Y}\|^{\frac{1}{2}}+C|\partial_z \tilde{Y}|\|Y\||\hat{Y}|^{\frac{1}{2}}\|\hat{Y}\|^{\frac{1}{2}}.
\end{eqnarray}
\end{lemma}

For the pressure term in (\ref{eq5-1}), define
\begin{eqnarray*}
G(Y):=P_H \left(                 
  \begin{array}{c}   
   -\int^z_{-h}\partial_x T dz'\\  
    0\\  
  \end{array}
\right), \quad {\rm{for}} \ Y=(v,T)\in V.
\end{eqnarray*}

Using the above functionals, we obtain
\begin{eqnarray}\label{equ-7}
\left\{
  \begin{array}{ll}
    dY(t)+AY(t)dt+B(Y(t),Y(t))dt+G(Y(t))dt=\psi(t,Y(t)) dW(t), \\
    Y(0)=Y_0,
  \end{array}
\right.
\end{eqnarray}
where
\begin{equation}\notag
W=\left(                 
  \begin{array}{c}   
    W_1\\  
    W_2 \\  
  \end{array}
\right) ,\quad
\psi(t,Y(t))
=\left(                 
  \begin{array}{cc}   
   \psi_1(t,Y(t)) & 0 \\  
   0 & \psi_2(t,Y(t)) \\  
  \end{array}
\right).
\end{equation}

\subsection{Definition of Strong Solution}
For the strong solution of (\ref{equ-7}), we shall fix a single stochastic basis $\mathcal{T}:=(\Omega, \mathcal{F}, \{\mathcal{F}_t\}_{t\geq 0}, P, W)$ with the expectation $\mathbb{E}$. Here,
\[
W=\left(                 
  \begin{array}{c}   
    W_1\\  
    W_2 \\  
  \end{array}
\right)
 \]
 is a cylindrical Brownian motion with the form $W(t,\omega)=\sum_{i\geq1}r_iw_i(t,\omega)$, where $\{r_i\}_{i\geq 1}$ is a complete orthonormal basis of a Hilbert space
$U=U_1\times U_2$, $U_1$ and $U_2$ are separable Hilbert spaces,
 $\{w_i\}_{i\geq1}$ is a sequence of independent one-dimensional standard Brownian motions on $(\Omega, \mathcal{F}, \{\mathcal{F}_t\}_{t\geq 0}, P)$.

In order to obtain the global well-posedness and moderate deviations of (\ref{equ-7}), we introduce the following Hypothesises:
\begin{description}
  \item[\textbf{Hypothesis A}] $\psi: [0,T]\times H\rightarrow \mathcal{L}_2({U}; H)$ satisfies that there exists a constant $K$ such that
\begin{description}
  \item[(A.1)] $\|\psi(t, \phi)\|^2_{\mathcal{L}_2({U}; H)}\leq K(1+|\phi|^2), \quad \forall t\in[0,T],\ \phi\in H$;
  \item[(A.2)] $\|\psi(t, \phi_1)-\psi(t,\phi_2)\|^2_{\mathcal{L}_2({U}; H)}\leq K|\phi_1-\phi_2|^2, \quad \forall t\in[0,T],\ \phi_1, \phi_2 \in H.$
\end{description}
\end{description}
\begin{description}
  \item[\textbf{Hypothesis B}] For the same constant $K$, we suppose that
 $\|\partial_z \psi(t, \phi)\|^2_{\mathcal{L}_2({U}; H)}\leq K(1+|\partial_z\phi|^2), \quad \forall t\in[0,T],\ \partial_z\phi\in H$.
\end{description}

Now, we give the definition of strong solution in probability  to (\ref{equ-7}).
\begin{dfn}\label{dfn-3}
Let $\mathcal{T}=(\Omega, \mathcal{F}, \{\mathcal{F}_t\}_{t\geq 0}, P, W)$ be a fixed stochastic basis and the initial condition $Y_0\in \mathcal{H}$.
An $\mathcal{F}_t-$predictable stochasitc process  $Y(t, \omega)$ is called a strong solution of (\ref{equ-7}) on [0,T] with the initial value $Y_0\in \mathcal{H}$ if for $P-$a.s. $\omega\in \Omega$,
\[
Y(\cdot)\in  C([0,T];H)\bigcap  L^2([0,T];V),\ \ \forall T>0,
\]
and for every $t\in [0, T]$,
\[
 (Y(t),\phi)-(Y_0,\phi)+\int^{t}_{0}\Big[\langle Y(s), A \phi\rangle+\langle B(Y,Y), \phi\rangle +( G(Y), \phi)\Big]ds=\int^{t}_{0}(\psi(s,Y(s)) dW(s),\phi), \quad P-a.s.
 \]
for all $\phi\in D(A)$.
\end{dfn}
To study the long-time behavior of the system (\ref{equ-7}), some kinds of $V-$norm estimates are needed, but it is difficult for our model. Fortunately, taking advantage of the special geometry structure of (\ref{equ-7}), we only need the moment estimate of $\partial_{z}v$ in $H$, which has been obtained during the proof of global well-posedness of the stochastic system (see \cite{G-S, GHZ08}).
\begin{thm}\label{thm-3}
Let the initial value $Y_0\in \mathcal{H}$. Assume that \textbf{Hypothesis A} and \textbf{Hypothesis B} hold, then there exists a unique global solution $Y$ of (\ref{equ-7}) in the sense of Definition  \ref{dfn-3} with $Y(0)=Y_0\in \mathcal{H}$. Furthermore, there exists a constant $C(K,T)$ such that
\begin{eqnarray}\label{ee-5}
\mathbb{E}\Big(\sup_{0\leq s\leq T}|Y(s)|^p+\int^T_0|Y(s)|^{p-2}\|Y(s)\|^2ds\Big)\leq C(K,T),
\end{eqnarray}
and
\begin{eqnarray}\label{ee-6}
\mathbb{E}\Big(\sup_{0\leq s\leq T}|\partial_zY(s)|^4+\int^T_0|\partial_zY(s)|^2\|\partial_zY(s)\|^2ds\Big)\leq C(K,T).
\end{eqnarray}
\end{thm}
Consider
\begin{eqnarray}\label{e-1}
\left\{
  \begin{array}{ll}
    dY^\varepsilon(t)+AY^\varepsilon(t)dt+B(Y^\varepsilon(t),Y^\varepsilon(t))dt+G(Y^\varepsilon(t))dt=\sqrt{\varepsilon}\psi(t,Y^\varepsilon(t)) dW(t), \\
    Y^\varepsilon(0)=Y_0\in \mathcal{H}.
  \end{array}
\right.
\end{eqnarray}

As the parameter $\varepsilon$ tends to 0, the solution $Y^\varepsilon$ of (\ref{e-1}) will tend to the solution of the following SPDE
\begin{eqnarray}\label{e-2}
\left\{
  \begin{array}{ll}
    dY^0(t)+AY^0(t)dt+B(Y^0(t),Y^0(t))dt+G(Y^0(t))dt=0, \\
    Y^0(0)=Y_0\in \mathcal{H}.
  \end{array}
\right.
\end{eqnarray}
As stated in the introduction, we will investigate the asymptotic behavior of the trajectory,
\begin{eqnarray}\label{eq-68}
Z^\varepsilon(t)=\frac{1}{\sqrt{\varepsilon}\lambda(\varepsilon)}(Y^\varepsilon-Y^0)(t), \ t\in[0,T],
\end{eqnarray}
where $\lambda(\varepsilon)$ is equal to 1  or satisfies (\ref{e-43}).
\section{Central limit theorem}\label{CLT}
In this part, we will estimate the central limit theorem, i.e. $\lambda(\varepsilon)=1$ in (\ref{eq-68}).

Let $Y^\varepsilon$ be the unique solution of (\ref{e-1}) in $L^2(\Omega; C([0,T];H))\cap L^2(\Omega; L^2([0,T]; V))$ and $Y^0$ be the unique solution of (\ref{e-2}). Taking an similar argument as in the proof of Theorem 3.1, we can also have
\begin{lemma}\label{lemma-2}
Let $Y_0\in \mathcal{H}$. Assume that \textbf{Hypothesis A} and \textbf{Hypothesis B} hold, then there exists $\varepsilon_0>0$ such that
\begin{eqnarray}\label{ee-7}
\sup_{\varepsilon \in [0, \varepsilon_0]}\mathbb{E}\Big(\sup_{0\leq s\leq T}|Y^\varepsilon(s)|^p+\int^T_0|Y^\varepsilon(s)|^{p-2}\|Y^\varepsilon(s)\|^2ds\Big)\leq C(K,T),
\end{eqnarray}
and
\begin{eqnarray}\label{ee-8}
\sup_{\varepsilon \in [0, \varepsilon_0]}\mathbb{E}\Big(\sup_{0\leq s\leq T}|\partial_zY^\varepsilon(s)|^4+\int^T_0|\partial_zY^\varepsilon(s)|^2\|\partial_zY^\varepsilon(s)\|^2ds\Big)\leq C(K,T).
\end{eqnarray}
In particular, it holds that
\begin{eqnarray}\label{ee-9}
\sup_{0\leq s\leq T}|Y^0(s)|^p+\int^T_0|Y^0(s)|^{p-2}\|Y^0(s)\|^2ds\leq C(K,T),
\end{eqnarray}
and
\begin{eqnarray}\label{ee-10}
\sup_{0\leq s\leq T}|\partial_zY^0(s)|^4+\int^T_0|\partial_zY^0(s)|^2\|\partial_zY^0(s)\|^2ds\leq C(K,T).
\end{eqnarray}
where $C(K,T)$ is a constant that depends only on $K$ and $T$.
\end{lemma}

Now, we explore the convergence of $Y^\varepsilon$ as $\varepsilon\rightarrow  0$.
\begin{prp}\label{prp-3}
Let $Y_0\in \mathcal{H}$. Assume that \textbf{Hypothesis A} and \textbf{Hypothesis B} hold, then there exists a constant $\varepsilon_0>0$ such that for all $0<\varepsilon\leq \varepsilon_0$,
\begin{eqnarray}
\mathbb{E}\left(\sup_{0\leq t\leq T}|Y^\varepsilon(t)-Y^0(t)|^2+\int^T_0\|Y^\varepsilon(s)-Y^0(s)\|^2ds\right)\leq \varepsilon C(K,T).
 \end{eqnarray}
\end{prp}
\begin{flushleft}
\textbf{Proof.} \quad Define $\tau_N=\inf\{t: |Y^\varepsilon(t)|^2+\int^t_0\|Y^\varepsilon(s)\|^2ds>N\}$.
Set $X^\varepsilon=Y^\varepsilon-Y^0$, we deduce from (\ref{e-1}) and (\ref{e-2}) that
\begin{eqnarray}\label{ee-12}
\left\{
  \begin{array}{ll}
    dX^\varepsilon(t)+AX^\varepsilon(t)dt+\Big(B(Y^\varepsilon(t),Y^\varepsilon(t))-B(Y^0,Y^0)\Big)dt+G(X^\varepsilon(t))dt=\sqrt{\varepsilon}\psi(t,Y^\varepsilon(t)) dW(t), \\
    X^\varepsilon(0)=0.
  \end{array}
\right.
\end{eqnarray}
Applying It\^{o} formula to $|X^\varepsilon|^2$ and by (\ref{ee-2}), we have
\begin{eqnarray}\notag
&&|X^\varepsilon(t\wedge \tau_N)|^2+2\int^{t\wedge \tau_N}_0\|X^\varepsilon(s)\|^2ds\\ \notag
&=&-2\int^{t\wedge \tau_N}_0\langle B(Y^\varepsilon, X^\varepsilon)+B(X^\varepsilon,Y^0), X^\varepsilon\rangle ds-2\int^{t\wedge \tau_N}_0(G(X^\varepsilon), X^\varepsilon)ds\\ \notag
&&+2\sqrt{\varepsilon}\int^{t\wedge \tau_N}_0(\psi(s,Y^\varepsilon)dW(s), X^\varepsilon)+\varepsilon\int^{t\wedge \tau_N}_0\|\psi(s,Y^\varepsilon)\|^2_{\mathcal{L}_2(U;H)}ds\\ \notag
&\leq& 2\int^{t\wedge \tau_N}_0|\langle B(X^\varepsilon,Y^0), X^\varepsilon\rangle | ds+2\int^{t\wedge \tau_N}_0|(G(X^\varepsilon), X^\varepsilon)|ds\\  \notag
&&+2\sqrt{\varepsilon}|\int^{t\wedge \tau_N}_0(\psi(s,Y^\varepsilon)dW(s), X^\varepsilon)|+\varepsilon\int^{t\wedge \tau_N}_0\|\psi(s,Y^\varepsilon)\|^2_{\mathcal{L}_2(U;H)}ds\\
\label{e-8}
&:=&I_1(t)+I_2(t)+I_3(t)+I_4(t).
\end{eqnarray}
Taking the supremum up to time $t\wedge \tau_N$ in (\ref{e-8}), and then taking the expectation, we have
\begin{eqnarray}\label{e-10}\notag
&&\mathbb{E}\left(\sup_{0\leq s\leq t\wedge \tau_N}\Big(|X^\varepsilon(s)|^2+2\int^{s}_0\|X^\varepsilon(l)\|^2dl\Big)\right)\\ \notag
&\leq& \mathbb{E}(I_1(t)+I_2(t)+\sup_{s\in[0,t\wedge \tau_N]}I_3(s)+I_4(t)).
\end{eqnarray}
By (\ref{ee-3}) and the Young's inequality, we have
\begin{eqnarray*}
\mathbb{E}I_1(t)
&\leq& C\mathbb{E}\int^{t\wedge \tau_N}_0(\|Y^0\||X^\varepsilon|\|X^\varepsilon\|+|\partial_z Y^0|\|X^\varepsilon\||X^\varepsilon|^{\frac{1}{2}}\|X^\varepsilon\|^{\frac{1}{2}})ds\\
&\leq& \frac{1}{4}\mathbb{E}\int^{t\wedge \tau_N}_0\|X^\varepsilon\|^2 ds+C\mathbb{E}\int^{t\wedge \tau_N}_0 (\|Y^0\|^2+|\partial_z Y^0|^4)|X^\varepsilon|^2ds.
\end{eqnarray*}
By the Cauchy-Schwarz inequality and the Young's inequality, we obtain
\begin{eqnarray*}
\mathbb{E}I_2(t)
&\leq& C\mathbb{E}\int^{t\wedge \tau_N}_0|X^\varepsilon|\|X^\varepsilon\|ds\\
&\leq& \frac{1}{4}\mathbb{E}\int^{t\wedge \tau_N}_0\|X^\varepsilon\|^2 ds+C\mathbb{E}\int^{t\wedge \tau_N}_0 |X^\varepsilon|^2ds.
\end{eqnarray*}
Moreover, by \textbf{Hypothesis $\mathbf{A}$}, we deduce that
\begin{eqnarray*}
\mathbb{E}I_4(t)\leq \varepsilon K\mathbb{E}\int^{t\wedge \tau_N}_0(1+|Y^\varepsilon|^2)ds.
\end{eqnarray*}
Applying the Burkholder-Davis-Gundy inequality and \textbf{Hypothesis $\mathbf{A}$}, we have
\begin{eqnarray*}
\mathbb{E}\sup_{s\in[0,t\wedge \tau_N]}I_3(s)
&\leq &2\sqrt{\varepsilon }\mathbb{E}\left(\int^{t\wedge \tau_N}_0|X^\varepsilon|^2\|\psi(s,Y^\varepsilon)\|^2_{\mathcal{L}_2(U;H)}ds\right)^{\frac{1}{2}}\\
&\leq& 4\sqrt{\varepsilon K}\mathbb{E}\left(\int^{t\wedge \tau_N}_0|X^\varepsilon|^2(1+|Y^\varepsilon|^2)ds\right)^{\frac{1}{2}}\\
&\leq& \frac{1}{2}\mathbb{E}\sup_{0\leq s\leq t\wedge \tau_N}|X^\varepsilon|^2 + 8\varepsilon K\mathbb{E}\int^{t\wedge \tau_N}_0(1+|Y^\varepsilon|^2)ds.
\end{eqnarray*}
Combing the above estimates and (\ref{ee-7}), we get
\begin{eqnarray}\label{e-13}\notag
&&\mathbb{E}\left(\sup_{0\leq s\leq t\wedge \tau_N}\Big(|X^\varepsilon(s)|^2+\int^{s}_0\|X^\varepsilon(l)\|^2dl\Big)\right)\\ \notag
&\leq& C\mathbb{E}\int^{t\wedge \tau_N}_0(\|Y^0\|^2+|\partial_z Y^0|^4+1)|X^\varepsilon|^2ds+C\varepsilon K\mathbb{E}\int^{t\wedge \tau_N}_0(1+|Y^\varepsilon|^2)ds\\
&\leq& C\mathbb{E}\int^{t\wedge \tau_N}_0(\|Y^0\|^2+|\partial_z Y^0|^4+1)|X^\varepsilon|^2ds+\varepsilon C(K,T).
\end{eqnarray}
By Gronwall inequality, (\ref{ee-9}) and (\ref{ee-10}), we have
\begin{eqnarray}\label{ee-11}\notag
&&\mathbb{E}\left(\sup_{0\leq s\leq t\wedge \tau_N}\Big(|X^\varepsilon(s)|^2+\int^{s}_0\|X^\varepsilon(l)\|^2dl\Big)\right)\\ \notag
&\leq& \varepsilon C(K,T)\cdot \exp\Big\{\int^{t\wedge \tau_N}_0(\|Y^0\|^2+|\partial_z Y^0|^4+1)ds\Big\}\\
&\leq& \varepsilon C(K,T).
\end{eqnarray}
Letting $N\rightarrow \infty$, we further obtain
\begin{eqnarray*}
\mathbb{E}\left(\sup_{0\leq t\leq T}|X^\varepsilon(t)|^2+\int^{T}_0\|X^\varepsilon(t)\|^2dt\right)
\leq \varepsilon C(K,T).
\end{eqnarray*}

\end{flushleft}
$\hfill\blacksquare$

Let $\tilde{V}^0=(\tilde{v}^0,\tilde{T}^0)$ be the solution of the following SPDE:
\begin{eqnarray}\label{e-21}
d\tilde{V}^0(t)+A\tilde{V}^0(t)dt+\Big(B(Y^0,\tilde{V}^0)+B( \tilde{V}^0(t),Y^0)\Big)dt+G(\tilde{V}^0(t))dt=\psi(t,Y^0(t))dW(t)
\end{eqnarray}
with $\tilde{V}^0(0)=0$. Under \textbf{Hypothesis A} and \textbf{Hypothesis B}, the global existence and uniqueness of the strong solution to (\ref{e-21}) can be proved easily by changing (\ref{e-21}) into a linear partial differential equation with random coefficients.  Taking a similar argument as in Lemma 4.1,  we have
\begin{eqnarray}\label{e-22}
\mathbb{E}\left(\sup_{0\leq t\leq T}|\tilde{V}^0|^p+\int^T_0|\tilde{V}^0|^{p-2}\|\tilde{V}^0(s)\|^2ds\right)\leq C(K,T).
\end{eqnarray}
and
\begin{eqnarray}\label{ee-14}
\mathbb{E}\left(\sup_{0\leq t\leq T}|\partial_z\tilde{V}^0|^4+\int^T_0|\partial_z\tilde{V}^0 |^2\|\partial_z\tilde{V}^0(s)\|^2ds\right)\leq C(K,T).
\end{eqnarray}

The following is the first result in this article: central limit theorem.
\begin{thm}(Central limit theorem)\label{thm-7}
Let $Y_0\in \mathcal{H}$. Assume \textbf{Hypothesis A} and \textbf{Hypothesis B} hold, then $\frac{Y^\varepsilon-Y^0}{\sqrt{\varepsilon}}$ converges to $\tilde{V}^0$ in the space $L^2(\Omega; C([0,T];H))\bigcap L^2(\Omega; L^2([0,T];V))$, that is
\begin{eqnarray}\label{e-23}
\lim_{\varepsilon\rightarrow 0}\mathbb{E}\left\{\sup_{0\leq t\leq T}\Big|\frac{Y^\varepsilon-Y^0}{\sqrt{\varepsilon}}-\tilde{V}^0\Big|^2+\int^T_0\Big\|\frac{Y^\varepsilon-Y^0}{\sqrt{\varepsilon}}-\tilde{V}^0(s)\Big\|^2ds\right\}=0.
\end{eqnarray}
\end{thm}
\begin{flushleft}
\textbf{Proof}. \quad  Let $\tilde{V}^\varepsilon(t):=\frac{Y^\varepsilon-Y^0}{\sqrt{\varepsilon}}=\frac{X^\varepsilon}{\sqrt{\varepsilon}}=(\tilde{v}^\varepsilon(t),\tilde{T}^\varepsilon(t))$, then $\tilde{V}^\varepsilon$ satisfies that
\begin{eqnarray}\label{e-24}
\left\{
  \begin{array}{ll}
    d\tilde{V}^\varepsilon(t)+A\tilde{V}^\varepsilon(t)dt+\Big(B(\tilde{V}^\varepsilon(t),Y^\varepsilon)+B(Y^0,\tilde{V}^\varepsilon(t))\Big)dt+G(\tilde{V}^\varepsilon(t))dt=\psi(t,Y^\varepsilon)dW(t), \\
 \tilde{V}^\varepsilon(0)=0.
  \end{array}
\right.
\end{eqnarray}
Define
\[
 \tilde{\rho}^\varepsilon(t)=\tilde{V}^\varepsilon(t)-\tilde{V}^0(t).
\]
From (\ref{e-21}) and (\ref{e-24}), we have
\begin{eqnarray*}
d\tilde{\rho}^\varepsilon(t)&+&A\tilde{\rho}^\varepsilon(t)dt+\Big(B(\tilde{V}^\varepsilon(t),Y^\varepsilon)-B(\tilde{V}^0(t),Y^0)\Big)dt
+B(Y^0,\tilde{\rho}^\varepsilon(t))dt+G(\tilde{\rho}^\varepsilon(t))dt\\
&&=(\psi(t,Y^\varepsilon)-\psi(t,Y^0))dW(t),
\end{eqnarray*}
Set
\[
\tau_N:=\inf\left\{t:|\tilde{\rho}^\varepsilon(t)|^2+\int^t_0\|\tilde{\rho}^\varepsilon(t)\|^2ds>N\right\}.
\]
 Applying It\^{o} formula to $|\tilde{\rho}^\varepsilon(t)|^2$ and by (\ref{ee-2}), we have
\begin{eqnarray*}
&&|\tilde{\rho}^\varepsilon(t\wedge{\tau_N})|^2+2\int^{t\wedge\tau_N}_0\|\tilde{\rho}^\varepsilon\|^2ds\\
 &=&-2\int^{t\wedge\tau_N}_0\langle B(\tilde{V}^\varepsilon, Y^\varepsilon-Y^0)+B(\tilde{\rho}^\varepsilon, Y^0), \tilde{\rho}^\varepsilon \rangle ds -2\int^{t\wedge\tau_N}_0(G(\tilde{\rho}^\varepsilon), \tilde{\rho}^\varepsilon)ds\\
&& \quad +2\int^{t\wedge\tau_N}_0((\psi(s,Y^\varepsilon)-\psi(s,Y^0))dW(s),\tilde{\rho}^\varepsilon)+\int^{t\wedge\tau_N}_0\|\psi(s,Y^\varepsilon)-\psi(s,Y^0)\|^2_{\mathcal{L}_2(U;H)}ds \\
&\leq &2\int^{t\wedge\tau_N}_0|\langle B(\tilde{V}^\varepsilon, Y^\varepsilon-Y^0), \tilde{\rho}^\varepsilon \rangle|ds
+2\int^{t\wedge\tau_N}_0|\langle B(\tilde{\rho}^\varepsilon, Y^0), \tilde{\rho}^\varepsilon \rangle| ds +2\int^{t\wedge\tau_N}_0|(G(\tilde{\rho}^\varepsilon), \tilde{\rho}^\varepsilon)|ds\\
&& \quad +2|\int^{t\wedge\tau_N}_0((\psi(s,Y^\varepsilon)-\psi(s,Y^0))dW(s),\tilde{\rho}^\varepsilon)|+\int^{t\wedge\tau_N}_0\|\psi(s,Y^\varepsilon)-\psi(s,Y^0)\|^2_{\mathcal{L}_2(U;H)}ds \\
&:=&J_1(t)+J_2(t)+J_3(t)+J_4(t)+J_5(t).
\end{eqnarray*}
Taking the supremum up to time $t\wedge \tau_N$ in the above equation, and taking the expectation, we obtain
\begin{eqnarray*}
&&\mathbb{E}\left(\sup_{0\leq s\leq t\wedge \tau_N}\Big(|\tilde{\rho}^\varepsilon(s)|^2+2\int^s_0\|\tilde{\rho}^\varepsilon(l)\|^2dl\Big)\right)\\
 &\leq & \mathbb{E}\Big(J_1(t)+J_2(t)+J_3(t)+\sup_{0\leq s\leq t}J_4(s)+J_5(t)\Big).
\end{eqnarray*}
By definition, we have $Y^\varepsilon-Y^0=\sqrt{\varepsilon}\tilde{V}^\varepsilon$. With the help of (\ref{ee-2}), (\ref{ee-3}), integration by parts and the Young's inequality, we have
\begin{eqnarray*}
 \mathbb{E}J_1(t)&=& 2\sqrt{\varepsilon}\mathbb{E}\int^{t\wedge \tau_N}_0|\langle B(\tilde{V}^\varepsilon,\tilde{V}^\varepsilon ), \tilde{\rho}^\varepsilon \rangle|ds\\
 &=& 2\sqrt{\varepsilon}\mathbb{E}\int^{t\wedge \tau_N}_0|\langle B(\tilde{V}^\varepsilon,\tilde{V}^\varepsilon ), \tilde{V}^0 \rangle|ds\\
 &=& 2\sqrt{\varepsilon}\mathbb{E}\int^{t\wedge \tau_N}_0|\langle B(\tilde{V}^\varepsilon,\tilde{V}^0), \tilde{V}^\varepsilon \rangle|ds\\
 &\leq& 2\sqrt{\varepsilon}\mathbb{E}\int^{t\wedge \tau_N}_0(\|\tilde{V}^0\||\tilde{V}^\varepsilon|\|\tilde{V}^\varepsilon\|+|\partial_z \tilde{V}^0|\|\tilde{V}^\varepsilon\||\tilde{V}^\varepsilon|^{\frac{1}{2}}\|\tilde{V}^\varepsilon\|^{\frac{1}{2}})ds\\
&\leq& C\sqrt{\varepsilon}\mathbb{E}\int^{t\wedge \tau_N}_0|\tilde{V}^\varepsilon|^2\|\tilde{V}^\varepsilon\|^2ds+C\sqrt{\varepsilon}\mathbb{E}\int^{t\wedge \tau_N}_0(\|\tilde{V}^0\|^2+\|\tilde{V}^\varepsilon\|^2+|\partial_z \tilde{V}^0|^4)ds\\
&\leq& C\sqrt{\varepsilon}\mathbb{E}\int^{t\wedge \tau_N}_0|\tilde{V}^\varepsilon|^2\|\tilde{V}^\varepsilon\|^2ds+\sqrt{\varepsilon} C(K,T).
\end{eqnarray*}
By (\ref{ee-3}), we get
\begin{eqnarray*}
 \mathbb{E}J_2(t)&\leq& 2C\mathbb{E}\int^{t\wedge \tau_N}_0(\|Y^0\||\tilde{\rho}^\varepsilon |\|\tilde{\rho}^\varepsilon \|+|\partial_z Y^0|\|\tilde{\rho}^\varepsilon \||\tilde{\rho}^\varepsilon |^{\frac{1}{2}}\|\tilde{\rho}^\varepsilon \|^{\frac{1}{2}})ds\\
 &\leq&\frac{1}{4}\mathbb{E}\int^{t\wedge \tau_N}_0 \|\tilde{\rho}^\varepsilon \|^2 ds+C\mathbb{E}\int^{t\wedge \tau_N}_0(\|Y^0\|^2+|\partial_z Y^0|^4)|\tilde{\rho}^\varepsilon|^2ds.
\end{eqnarray*}
By Cauchy-Schwarz inequality and the Young's inequality, we otain
\begin{eqnarray*}
 \mathbb{E}J_3(t)&\leq& \frac{1}{4}\mathbb{E}\int^{t\wedge \tau_N}_0 \|\tilde{\rho}^\varepsilon \|^2 ds+C\mathbb{E}\int^{t\wedge \tau_N}_0|\tilde{\rho}^\varepsilon|^2ds.
\end{eqnarray*}
Applying the Burkholder-Davis-Gundy inequality to $J_4(t)$ and by \textbf{Hypothesis A}, we obtain
\begin{eqnarray*}
 \mathbb{E}\sup_{0\leq s\leq t\wedge \tau_N}J_4(s)
 &\leq& 4C\mathbb{E}\left(\int^{t\wedge \tau_N}_0|\tilde{\rho}^\varepsilon(s)|^2\|\psi(s, Y^\varepsilon)-\psi(s, Y^0)\|^2_{\mathcal{L}_2(U;H)}ds\right)^{\frac{1}{2}}\\
 &\leq& 4C\sqrt{K}\mathbb{E}(\int^{t\wedge \tau_N}_0|\tilde{\rho}^\varepsilon(s)|^2|X^\varepsilon|^2ds)^{\frac{1}{2}}\\
 &\leq& \frac{1}{2}\mathbb{E}\sup_{0\leq s\leq t\wedge \tau_N} |\tilde{\rho}^\varepsilon(s)|^2+C(K)\mathbb{E}\int^{t\wedge \tau_N}_0|X^\varepsilon|^2ds.
\end{eqnarray*}
By \textbf{Hypothesis A}, we have
\begin{eqnarray*}
 \mathbb{E}J_5(t)\leq K\mathbb{E}\int^{t\wedge \tau_N}_0|X^\varepsilon|^2ds.
\end{eqnarray*}
Collecting the above estimates and by (\ref{ee-7}), (\ref{ee-9}), we have
\begin{eqnarray*}\notag
&&\mathbb{E}\left(\sup_{0\leq s\leq T\wedge \tau_N}\Big(|\tilde{\rho}^\varepsilon(s)|^2+\int^s_0\|\tilde{\rho}^\varepsilon(l)\|^2dl\Big)\right)\\ \notag
 &\leq &
 C\sqrt{\varepsilon}\mathbb{E}\int^{T\wedge \tau_N}_0| \tilde{V}^\varepsilon|^2\|\tilde{V}^\varepsilon\|^2ds
+C\sqrt{\varepsilon}\mathbb{E}\int^{T\wedge \tau_N}_0(1+\|Y^0\|^2+|\partial_z Y^0|^4)|\tilde{\rho}^\varepsilon|^2ds+
\sqrt{\varepsilon}C(K,T).
\end{eqnarray*}
Suppose that there exists a constant $\varepsilon_0>0$ such that
\begin{eqnarray}\label{e-27}
\sup_{0\leq \varepsilon \leq \varepsilon_0}\mathbb{E}\int^{T}_0| \tilde{V}^\varepsilon|^2\|\tilde{V}^\varepsilon\|^2ds<C(K,T),
\end{eqnarray}
which will be proved in the following Lemma \ref{lem-3}. Then,
\begin{eqnarray}\notag
&&\mathbb{E}\left(\sup_{0\leq s\leq T\wedge \tau_N}\Big(|\tilde{\rho}^\varepsilon(s)|^2+\int^s_0\|\tilde{\rho}^\varepsilon(l)\|^2dl\Big)\right)\\
\label{e-25}
 &\leq &C\sqrt{\varepsilon}\mathbb{E}\int^{T\wedge \tau_N}_0(1+\|Y^0\|^2+|\partial_z Y^0|^4)|\tilde{\rho}^\varepsilon|^2ds+
\sqrt{\varepsilon}C(K,T).
\end{eqnarray}
Applying Gronwall inequality to (\ref{e-25}) and by (\ref{ee-9}), (\ref{ee-10}), we obtain
\begin{eqnarray*}\notag
&&\mathbb{E}\left(\sup_{0\leq s\leq T\wedge \tau_N}\Big(|\tilde{\rho}^\varepsilon(s)|^2+\int^s_0\|\tilde{\rho}^\varepsilon(l)\|^2dl\Big)\right)\\
 &\leq &\varepsilon C(K,T)\cdot \exp\Big\{\int^T_0 (1+\|Y^0\|^2+|\partial_z Y^0|^4) ds\Big\}\\
 &\leq &\varepsilon C(K,T).
\end{eqnarray*}
Letting $N\rightarrow \infty$, we
conclude the result.
\end{flushleft}
$\hfill\blacksquare$

To finish the proof of Theorem \ref{thm-7}, it remains to prove (\ref{e-27}).
\begin{lemma}\label{lem-3}
Let $Y_0\in \mathcal{H}$. Assume \textbf{Hypothesis A} and \textbf{Hypothesis B} hold, then there exists a constant $\varepsilon_0>0$ such that
\begin{eqnarray*}
\sup_{0\leq \varepsilon \leq \varepsilon_0}\mathbb{E}\int^{T}_0| \tilde{V}^\varepsilon|^2\|\tilde{V}^\varepsilon\|^2ds<C(K,T).
\end{eqnarray*}
\end{lemma}
\begin{flushleft}
\textbf{Proof.} \quad Recall that
\begin{eqnarray}\label{ee-15}
d\tilde{V}^\varepsilon(t)+A\tilde{V}^\varepsilon(t)dt+\Big(B(\tilde{V}^\varepsilon(t),Y^\varepsilon)+B(Y^0,\tilde{V}^\varepsilon(t))\Big)dt+G(\tilde{V}^\varepsilon(t))dt=\psi(t,Y^\varepsilon)dW(t),
\end{eqnarray}
with $\tilde{V}^\varepsilon(0)=0$.
Applying It\^{o} formula to $|\tilde{V}^\varepsilon|^2$ and by (\ref{ee-2}), we have
\begin{eqnarray*}
d|\tilde{V}^\varepsilon(t)|^2
&=&-2\|\tilde{V}^\varepsilon\|^2dt-2\langle B(\tilde{V}^\varepsilon(t),Y^\varepsilon)+B(Y^0,\tilde{V}^\varepsilon(t)), \tilde{V}^\varepsilon\rangle dt\\
&&
-2\left(G(\tilde{V}^\varepsilon(t)),\tilde{V}^\varepsilon\right)dt
+2(\psi(t, Y^\varepsilon)dW(t),\tilde{V}^\varepsilon)+\|\psi(t, Y^\varepsilon)\|^2_{\mathcal{L}_2(U;H)}dt,
\end{eqnarray*}
then
\[
\langle|\tilde{V}^\varepsilon(\cdot)|^2\rangle_t=4\int^t_0\|\psi(s, Y^\varepsilon)\|^2_{\mathcal{L}_2(U;H)}|\tilde{V}^\varepsilon|^2ds.
\]
Applying It\^{o} formula to $|\tilde{V}^\varepsilon|^4$, we obtain
\begin{eqnarray*}
d|\tilde{V}^\varepsilon(t)|^4
&=&2|\tilde{V}^\varepsilon|^2d|\tilde{V}^\varepsilon|^2+d\langle|\tilde{V}^\varepsilon(\cdot)|^2\rangle_t\\
&=& 2|\tilde{V}^\varepsilon|^2\Big[-2\|\tilde{V}^\varepsilon\|^2dt-2\langle B(\tilde{V}^\varepsilon(t),Y^\varepsilon), \tilde{V}^\varepsilon \rangle dt-2\left(G(\tilde{V}^\varepsilon(t)),\tilde{V}^\varepsilon\right)dt
\\
&&
+2(\psi(t, Y^\varepsilon)dW(t),\tilde{V}^\varepsilon)+\|\psi(t, Y^\varepsilon)\|^2_{\mathcal{L}_2(U;H)}dt\Big]+4\|\psi(t, Y^\varepsilon)\|^2_{\mathcal{L}_2(U;H)}|\tilde{V}^\varepsilon(t)|^2dt.
\end{eqnarray*}
Define $\tau_N:=\inf\left\{t: |\tilde{V}^\varepsilon(t)|^4+\int^t_0 \|\tilde{V}^\varepsilon\|^2ds>N\right\}$, then
\begin{eqnarray*}
&&|\tilde{V}^\varepsilon({t\wedge \tau_N})|^4+4\int^{t\wedge \tau_N}_0|\tilde{V}^\varepsilon|^2\|\tilde{V}^\varepsilon\|^2ds\\
&\leq& 4\int^{t\wedge \tau_N}_0|\tilde{V}^\varepsilon|^2|\langle B(\tilde{V}^\varepsilon,Y^\varepsilon), \tilde{V}^\varepsilon \rangle| ds+4\int^{t\wedge \tau_N}_0|\tilde{V}^\varepsilon|^2|\left(G(\tilde{V}^\varepsilon),\tilde{V}^\varepsilon\right)|ds\\
&&+6\int^{t\wedge \tau_N}_0|\tilde{V}^\varepsilon|^2\|\psi(s, Y^\varepsilon)\|^2_{\mathcal{L}_2(U;H)}ds+4|\int^{t\wedge \tau_N}_0|\tilde{V}^\varepsilon|^2(\psi(s, Y^\varepsilon)dW(s),\tilde{V}^\varepsilon)|\\
&:=& K_1(t)+K_2(t)+K_3(t)+K_4(t).
\end{eqnarray*}
Taking the expectation, we obtain
\begin{eqnarray}\notag \label{e-35}
&&\mathbb{E}\left(\sup_{0\leq s\leq {t\wedge \tau_N}}\Big(|\tilde{V}^\varepsilon(s)|^4+4\int^{s\wedge \tau_N}_0|\tilde{V}^\varepsilon|^2\|\tilde{V}^\varepsilon\|^2dl\Big)\right)\\
&\leq& \mathbb{E}K_1(t)+\mathbb{E}K_2(t)+\mathbb{E}K_3(t)+\mathbb{E}\sup_{0\leq s\leq {t\wedge \tau_N}}|K_4(s)|.
\end{eqnarray}
Notice that
\begin{eqnarray}\label{eee-1}
Y^\varepsilon=Y^0+\sqrt{\varepsilon}\tilde{V}^\varepsilon,
\end{eqnarray}
then by (\ref{ee-2}), we get
\begin{eqnarray*}
 \mathbb{E}K_1(t)&=&4\mathbb{E}\int^{t\wedge \tau_N}_0|\tilde{V}^\varepsilon|^2|\langle B(\tilde{V}^\varepsilon,Y^0+\sqrt{\varepsilon}\tilde{V}^\varepsilon), \tilde{V}^\varepsilon\rangle|ds\\
 &=&4\mathbb{E}\int^{t\wedge \tau_N}_0|\tilde{V}^\varepsilon|^2|\langle B(\tilde{V}^\varepsilon,Y^0), \tilde{V}^\varepsilon\rangle|ds.
\end{eqnarray*}
By (\ref{ee-3}) and the Young's inequality, we deduce that
\begin{eqnarray}\notag
\mathbb{E}K_1(t)&\leq& 4\mathbb{E}\int^{t\wedge \tau_N}_0|\tilde{V}^\varepsilon|^2(\|\tilde{V}^\varepsilon\||\tilde{V}^\varepsilon| \|Y^{0}\|+|\partial_z Y^{0}|\|\tilde{V}^\varepsilon\|^{\frac{3}{2}} |\tilde{V}^\varepsilon|^{\frac{1}{2}})ds\\ \notag
&\leq& 2\mathbb{E}\int^{t\wedge \tau_N}_0|\tilde{V}^\varepsilon|^2\| \tilde{V}^\varepsilon\|^2ds
+C\mathbb{E}\int^{t\wedge \tau_N}_0(\|Y^0\|^2+|\partial_{z}Y^{0}|^{4})|\tilde{V}^\varepsilon|^4ds.
\end{eqnarray}

Using the Cauchy-Schwarz inequality and the Young's inequality, we obtain
\begin{eqnarray*}\notag
\mathbb{E}K_2(t)
&\leq&4\mathbb{E}\int^{t\wedge \tau_N}_0|\tilde{V}^\varepsilon|^3\|\tilde{V}^\varepsilon\|ds\\
\label{e-38}
&\leq& \frac{1}{4}\mathbb{E}\int^{t\wedge \tau_N}_0|\tilde{V}^\varepsilon|^2\|\tilde{V}^\varepsilon\|^2ds+C\mathbb{E}\int^{t\wedge \tau_N}_0|\tilde{V}^\varepsilon|^4ds.
\end{eqnarray*}

By \textbf{Hypothesis A}, (\ref{ee-7}), (\ref{ee-9}) and (\ref{eee-1}), we have
 \begin{eqnarray}\notag
\mathbb{E}K_3(t)
&\leq&6K\mathbb{E}\int^{t\wedge \tau_N}_0(1+|Y^\varepsilon|^2)|\tilde{V}^\varepsilon|^2ds\\ \notag
&\leq& 6K\mathbb{E}\int^{t\wedge \tau_N}_0|\tilde{V}^\varepsilon|^2ds+12K\mathbb{E}\int^{t\wedge \tau_N}_0(|Y^0|^2+\varepsilon|\tilde{V}^\varepsilon|^2)|\tilde{V}^\varepsilon|^2ds\\
\label{e-39}
&\leq& 12K\varepsilon\mathbb{E}\int^{t\wedge \tau_N}_0|\tilde{V}^\varepsilon|^4ds+C(K,T).
\end{eqnarray}
Applying the Burkholder-Davis-Gundy inequality, the Young's inequality, (\ref{ee-9}) and (\ref{eee-1}), we deduce that
\begin{eqnarray}\notag
\mathbb{E}\sup_{0\leq s\leq {t\wedge \tau_N}}|K_4(s)|
&\leq& C\mathbb{E}\left(\int^{t\wedge \tau_N}_0|\tilde{V}^\varepsilon|^4\|\psi(s, Y^\varepsilon)\|^2_{\mathcal{L}_2(U;H)}|\tilde{V}^\varepsilon|^2ds\right)^{\frac{1}{2}}\\ \notag
&\leq& C\mathbb{E}\left[\sup_{0\leq s\leq {t\wedge \tau_N}}|\tilde{V}^\varepsilon|^2 \cdot \left(\int^{t\wedge \tau_N}_0\|\psi(s, Y^\varepsilon)\|^2_{\mathcal{L}_2(U;V)}|\tilde{V}^\varepsilon|^2ds\right)^{\frac{1}{2}}\right]\\ \notag
&\leq& \frac{1}{4}\mathbb{E}\sup_{0\leq s\leq t\wedge \tau_N}|\tilde{V}^\varepsilon(s)|^4+CK\mathbb{E}\int^{t\wedge \tau_N}_0(1+|Y^\varepsilon|^2)|\tilde{V}^\varepsilon|^2ds\\ \notag
&\leq&\frac{1}{4}\mathbb{E}\sup_{0\leq s\leq t\wedge \tau_N}|\tilde{V}^\varepsilon(s)|^4+CK\mathbb{E}\int^{t\wedge \tau_N}_0(1+|Y^0|^2+\varepsilon |\tilde{V}^\varepsilon|^2)|\tilde{V}^\varepsilon|^2ds\\
\label{e-40}
&\leq&\frac{1}{4}\mathbb{E}\sup_{0\leq s\leq t\wedge \tau_N}|\tilde{V}^\varepsilon(s)|^4+\varepsilon CK\mathbb{E}\int^{t\wedge \tau_N}_0|\tilde{V}^\varepsilon|^4 ds+C(K,T).
\end{eqnarray}
As a result of (\ref{e-35})-(\ref{e-40}), we have
\begin{eqnarray}\notag
&&\mathbb{E}\left[\sup_{0\leq s\leq {t\wedge \tau_N}}\Big(\frac{1}{2}|\tilde{V}^\varepsilon(s)|^4+\int^{s\wedge \tau_N}_0|\tilde{V}^\varepsilon|^2\|\tilde{V}^\varepsilon\|^2dl\Big)\right]\\
\label{e-41}
&\leq& C\mathbb{E}\int^{t\wedge \tau_N}_0(\|Y^0\|^2+|\partial_z Y^0|^4+CK\varepsilon+1)|\tilde{V}^\varepsilon|^4ds+C(K,T).
\end{eqnarray}
By Gronwall inequality, (\ref{ee-9}) and (\ref{ee-10}), we have
\begin{eqnarray*}\notag
&&\mathbb{E}\left[\sup_{0\leq s\leq {t\wedge \tau_N}}\Big(|\tilde{V}^\varepsilon(s)|^4+2\int^{s\wedge \tau_N}_0|\tilde{V}^\varepsilon|^2\|\tilde{V}^\varepsilon\|^2dl\Big)\right]\\
&\leq& C(K,T)\cdot \exp\Big\{\int^{T}_0(\|Y^0\|^2+|\partial_z Y^0|^4+CK\varepsilon+1)ds\Big\}\\
&\leq& C(K,T)(1+\varepsilon).
\end{eqnarray*}
Let $N\rightarrow \infty$, for a constant $\varepsilon_0>0$,  we get
\begin{eqnarray*}
\sup_{0\leq \varepsilon \leq \varepsilon_0}\mathbb{E}\left[\sup_{0\leq s\leq T}|\tilde{V}^\varepsilon(s)|^4+2\int^T_0|\tilde{V}^\varepsilon|^2\|\tilde{V}^\varepsilon\|^2ds\right]\leq C(K,T),
\end{eqnarray*}
which implies the result.
\end{flushleft}
$\hfill\blacksquare$
\section{Moderate deviation principle}\label{MDP}
In this part, we are concerned with the moderate deviation principle of $Y^\varepsilon$. As stated in introduction, we need to
prove $\frac{Y^\varepsilon-Y^0}{\sqrt{\varepsilon}\lambda(\varepsilon)}$ satisfies a large deviation principle on $C([0,T]; H)\cap L^2([0,T];V)$ with $\lambda(\varepsilon)$ satisfying (\ref{e-43}).
From now on, we assume (\ref{e-43}) holds.
\subsection{The weak convergence approach}
Let $Z^\varepsilon =\frac{Y^\varepsilon-Y^0}{\sqrt{\varepsilon}\lambda(\varepsilon)}$, we will use the weak convergence approach introduced by Budhiraja and Dupuis in \cite{BD} to verify $Z^\varepsilon$ satisfies a large deviation principle. Firstly recall some standard definitions and results from the large deviation theory (see \cite{DZ}).

Suppose $\{Z^\varepsilon\}$ be a family of random variables defined on a probability space $(\Omega, \mathcal{F}, P)$ taking values in some Polish space $\mathcal{E}$.

\begin{dfn}
(Rate Function) A function $I: \mathcal{E}\rightarrow [0,\infty]$ is called a rate function if $I$ is lower semicontinuous. A rate function $I$ is called a good rate function if the level set $\{x\in \mathcal{E}: I(x)\leq M\}$ is compact for each $M<\infty$.
\end{dfn}
\begin{dfn}
\begin{description}
  \item[(i)] (Large deviation principle) The sequence $\{Z^{\varepsilon}\}$ is said to satisfy the large deviation principle with rate function $I$ if for each Borel subset $A$ of $\mathcal{E}$
      \[
      -\inf_{x\in A^o}I(x)\leq \lim \inf_{\varepsilon\rightarrow 0}\varepsilon \log P(Z^{\varepsilon}\in A)\leq \lim \sup_{\varepsilon\rightarrow 0}\varepsilon \log P(Z^{\varepsilon}\in A)\leq -\inf_{x\in \bar{A}}I(x),
      \]
      where $A^o$ and $\bar{A}$ denote the interior and closure of $A$ in $\mathcal{E}$, respectively.
  \item[(ii)] (Laplace principle) The sequence $\{Z^{\varepsilon}\}$ is said to satisfy the Laplace principle with rate function $I$ if for each bounded continuous real-valued function $f$ defined on $\mathcal{E}$
      \[
      \lim_{\varepsilon\rightarrow 0}\varepsilon \log E\Big\{\exp[-\frac{1}{\varepsilon}f(Z^{\varepsilon})]\Big\}=-\inf_{x\in \mathcal{E}}\{f(x)+I(x)\}.
      \]
\end{description}
\end{dfn}
It's well-known that the large deviation principle and the Laplace principle are equivalent if $\mathcal{E}$ is a Polish space and the rate function is good (see \cite{DZ}).

Suppose $W(t)$ is a cylindrical Wiener process on a Hilbert space $U$ defined on a probability space $(\Omega, \mathcal{F},\{\mathcal{F}_t\}_{t\in [0,T]}, P )$ ( the paths of $W$ take values in $C([0,T];\mathcal{U})$, where $\mathcal{U}$ is another Hilbert space such that the embedding $U\subset \mathcal{U}$ is Hilbert-Schmidt).
Now we define
\begin{eqnarray*}
&\mathcal{A}=\{\phi: \phi\ {\rm{is\ a \ U-valued  \ \{\mathcal{F}_t\}-predictable\ process}}\ s.t.\ \int^T_0 |\phi(s)|^2_Uds<\infty\ a.s.\};\\
&T_M=\{ h\in L^2([0,T];U): \int^T_0 |h(s)|^2_Uds\leq M\};\\
&\mathcal{A}_M=\{\phi\in \mathcal{A}: \phi(\omega)\in T_M,\ P\text{-}a.s.\}.
\end{eqnarray*}
Here, we use the weak topology on $L^2([0,T];U)$ under which $T_M$ is a compact space.

Suppose $\mathcal{G}^{\varepsilon}: C([0,T];\mathcal{U})\rightarrow \mathcal{E}$ is a measurable map and $Z^{\varepsilon}=\mathcal{G}^{\varepsilon}(W)$. Now, we list the following sufficient conditions for the Laplace principle (equivalently, large deviation principle) of $Z^{\varepsilon}$.
\begin{description}
  \item[\textbf{Hypothesis H1} ] There exists a measurable map $\mathcal{G}^0: C([0,T];\mathcal{U})\rightarrow \mathcal{E}$ satisfying
\end{description}
\begin{description}
  \item[(i)] For every $M<\infty$, let $\{h^{\varepsilon}: \varepsilon>0\}$ $\subset \mathcal{A}_M$. If $h^{\varepsilon}$ converges to $h$ as $T_M-$valued random elements in distribution, then $\mathcal{G}^{\varepsilon}(W(\cdot)+\lambda(\varepsilon)\int^{\cdot}_{0}h^\varepsilon(s)ds)$ converges in distribution to $\mathcal{G}^0(\int^{\cdot}_{0}h(s)ds)$.
  \item[(ii)] For every $M<\infty$, $K_M=\{\mathcal{G}^0(\int^{\cdot}_{0}h(s)ds): h\in T_M\}$ is a compact subset of $\mathcal{E}$.
\end{description}
The following result is due to Budhiraja et al. in \cite{BD}.

\begin{thm}\label{thm-2}
If $\mathcal{G}^{0}$ satisfies \textbf{Hypothesis H1}, then $Z^{\varepsilon}$ satisfies a large deviation principle on $\mathcal{E}$ with the good rate function $I$ given by
\begin{eqnarray}\label{eq-5}
I(f)=\inf_{\{h\in L^2([0,T];U): f= \mathcal{G}^0(\int^{\cdot}_{0}h(s)ds)\}}\Big\{\frac{1}{2}\int^T_0|h(s)|^2_{U}ds\Big\},\ \ \forall f\in\mathcal{E}.
\end{eqnarray}
By convention, $I(\emptyset)=\infty$.
\end{thm}
From (\ref{e-24}), $Z^\varepsilon$ satisfies the following SPDE
\begin{eqnarray}\label{e-45}
\left\{
  \begin{array}{ll}
     dZ^\varepsilon(t)+AZ^\varepsilon(t)dt+B(Z^\varepsilon(t),Y^0+\sqrt{\varepsilon}\lambda(\varepsilon)Z^\varepsilon)dt+B(Y^0 ,Z^\varepsilon(t))dt\\
     \quad\quad \quad\quad\quad\quad\quad\quad+G(Z^\varepsilon(t))dt=\lambda^{-1}(\varepsilon)\psi(t,Y^0+\sqrt{\varepsilon}\lambda(\varepsilon){Z}^\varepsilon)dW(t), \\
     Z^\varepsilon(0)=0.
  \end{array}
\right.
\end{eqnarray}

The following is our main result in this part.
\begin{thm}\label{thm-8}
Let $Y_0\in V$. Assume \textbf{Hypothesis A} and \textbf{Hypothesis B} hold. Then $Z^{\varepsilon}$ satisfies a large deviation principle on $C([0,T];H)\cap L^2([0,T];V)$ with the good rate function $I$ defined by (\ref{eq-5}) with respect to the uniform convergence.
\end{thm}
\subsection{Priori estimates}\label{Sec 4}

Under \textbf{Hypothesis A} and \textbf{Hypothesis B}, by Theorem \ref{thm-3},
there exists a pathwise unique strong solution $Z^\varepsilon$ of (\ref{e-45}) in $\Re:= C([0,T];H)\cap L^2([0,T]; V) $, the norm of $\Re$ is
\[
|Y|^2_{\Re}:=\sup_{0\leq t\leq T}|Y(t)|^2+\int^T_0\|Y(t)\|^2dt.
\]

Therefore, there exist Borel-measurable functions
\begin{eqnarray}\label{e-46}
\mathcal{G}^{\varepsilon}: C([0,T];\mathcal{U})\rightarrow \Re\text{  such that }Z^{\varepsilon}(\cdot)=\mathcal{G}^{\varepsilon}(W(\cdot)).
\end{eqnarray}




For $h\in L^2([0,T];U)$, consider the following skeleton equation of (\ref{e-45})
\begin{eqnarray}\label{equ-9}
\left\{
  \begin{array}{ll}
dR^h(t)+AR^h(t)dt+B(R^h,Y^0)dt+B(Y^0, R^h(t))dt+G(R^h(t))dt=\psi(t,Y^0(t))h(t)dt,\\
R^h(0)=0.
  \end{array}
\right.
\end{eqnarray}
\subsubsection{Global Well-posedness }
\begin{thm}\label{thm-4}
Assume \textbf{Hypothesis A} holds and the initial data $Y_0=(v_0,T_0)\in V$, $h\in T_M $, for some $M>0$, then for any $T>0$, (\ref{equ-9}) has a unique strong solution $R_h\in C([0,T];V)\bigcap L^2([0,T];D(A))$. Moreover,  there exists a constant $C(M,K,T)$ such that
\begin{eqnarray}\label{e-44}
\sup_{h\in T_M}\Big\{\sup_{0\leq t\leq T}\|R^h(t)\|^2+\int^T_0\|R^h(t)\|^2_{D(A)}dt \Big\}\leq C(M,K,T).
\end{eqnarray}
\end{thm}
\begin{flushleft}
\textbf{Proof.}\quad
Since the equation (\ref{equ-9}) is a linear equation, the proof of global well-posedness is standard. Therefore, we only prove (\ref{e-44}).
When $Y_0\in V$, referring to \cite{C-T-1}, we have
\begin{eqnarray}\label{e-80}
\sup_{0\leq t\leq T}\|Y^0(t)\|^2+\int^T_0\|Y^0(t)\|^2_{D(A)}dt \leq C(T).
\end{eqnarray}
Based on (\ref{e-80}) and similar to the process of proving the global well-posedness of the skeleton equation in \cite{RR}, we deduce that (\ref{e-44}) holds.

\end{flushleft}
$\hfill\blacksquare$

Now, we can define $\mathcal{G}^0: C([0,T];\mathcal{U})\rightarrow \Re$ by
\begin{eqnarray}
\mathcal{G}^0(\tilde{h})=\left\{
                 \begin{array}{ll}
                   R^h, & {\rm{if}}\  \tilde{h}=\int^{\cdot}_{0}h(s)ds \quad {\rm{for\ some\ }}  h\in L^2([0,T]; U), \\
                   0, & {\rm{otherwise}}.
                 \end{array}
               \right.
\end{eqnarray}

\subsubsection{Compactness of $R^h$}
In order to prove compactness of $R^h$, as in \cite{FG95}, we introduce the following space. Let $K$ be a separable Hilbert space. Given $p>1, \alpha\in (0,1)$, let $W^{\alpha,p}([0,T]; K)$ be the Sobolev space of all $u\in L^p([0,T];K)$ such that
\[
\int^T_0\int^T_0\frac{\|u(t)-u(s)\|_K^{ p}}{|t-s|^{1+\alpha p}}dtds< \infty,
\]
endowed with the norm
\[
\|u\|^p_{W^{\alpha,p}([0,T]; K)}=\int^T_0\|u(t)\|_K^pdt+\int^T_0\int^T_0\frac{\|u(t)-u(s)\|_K^{ p}}{|t-s|^{1+\alpha p}}dtds.
\]
The following result can be found in \cite{FG95}.
\begin{lemma}\label{lem-1}
Let $B_0\subset B\subset B_1$ be Banach spaces, $B_0$ and $B_1$ reflexive, with compact embedding of $B_0$ in $B$. Let $p\in (1, \infty)$ and $\alpha \in (0, 1)$ be given. Let $\Lambda$ be the space
\[
\Lambda= L^p([0, T]; B_0)\cap W^{\alpha, p}([0,T]; B_1),
\]
endowed with the natural norm. Then the embedding of $\Lambda$ in $L^p([0,T];B)$ is compact.
\end{lemma}
The following is the result we obtain in this part.
\begin{prp}\label{prp-6-1}
Let $Y_0 \in V$, under \textbf{Hypothesis A} and \textbf{Hypothesis B}, $R^h$ is compact in $L^2([0,T];V).$
\end{prp}
\begin{flushleft}
\textbf{Proof.} \quad
Define
\[
F(R^h,Y^0):=AR^h+B(R^h,Y^0)+B(Y^0,R^h)+G(R^h).
\]
From (\ref{equ-9}), we have
\begin{eqnarray*}
R^h(t)&=&-\int^t_0F(R^h(s),Y^0(s))ds
+\int^t_0\psi(s,Y^0)h(s)ds\\
&:=&I_1(t)+I_2(t).
\end{eqnarray*}
With the aid of the H\"{o}lder inequality and the Cauchy-Schwarz inequality, we have
\begin{eqnarray*}
\|F(R^h,Y^0)\|_{V'}\leq C\|R^h\|+C\|R^h\|\|Y^0\|+C\|R^h\||\partial_z Y^0|+C\|Y^0\||\partial_z R^h|+C|R^h|.
\end{eqnarray*}
By the H\"{o}lder inequality, we obtain
\begin{eqnarray*}
\|I_1(t)-I_1(s)\|^2_{V'}
&=&\|\int^t_s F(R^h,Y^0)dl\|^2_{V'}\\
&\leq& C(t-s)\int^t_s\| F(R^h,Y^0)\|^2_{V'}dl\\
&\leq& C(t-s)^2\Big[\sup_{t\in[0,T]}\|R^h(t)\|^2+\sup_{t\in[0,T]}\|R^h(t)\|^2\sup_{t\in[0,T]}\|Y^0(t)\|^2
\\
&&\ +\sup_{t\in[0,T]}|\partial_z Y^0|^2\sup_{t\in[0,T]}\|R^h(t)\|^2+\sup_{t\in[0,T]}|\partial_z R^h|^2\sup_{t\in[0,T]}\|Y^0(t)\|^2+\sup_{t\in[0,T]}|R^h(t)|^2\Big].
\end{eqnarray*}
From the definition of $W^{r,2}([0,T]; V')$, for $r\in (0,\frac{1}{2})$, we have
\begin{eqnarray*}
\|I_1\|^2_{W^{r,2}([0,T];V')}
&=&\int^T_0\|I_1(t)\|^2_{V'}dt+\int^T_0\int^T_0\frac{\|I_1(t)-I_1(s)\|^2_{V'}}{|t-s|^{1+2r}}dsdt\\
&\leq&
C(r, T)\Big[\sup_{t\in[0,T]}\|R^h(t)\|^2+\sup_{t\in[0,T]}\|R^h(t)\|^2\sup_{t\in[0,T]}\|Y^0(t)\|^2
\\
&&\ +\sup_{t\in[0,T]}|\partial_z Y^0|^2\sup_{t\in[0,T]}\|R^h(t)\|^2+\sup_{t\in[0,T]}|\partial_z R^h|^2\sup_{t\in[0,T]}\|Y^0(t)\|^2+\sup_{t\in[0,T]}|R^h(t)|^2\Big].
\end{eqnarray*}
Thus, by (\ref{e-44}) and (\ref{e-80}), we have
\begin{eqnarray*}
\mathbb{E}\|I_1\|_{W^{r,2}([0,T];V')}
\leq C_{1}(r,M,K,T).
\end{eqnarray*}
By \textbf{Hypothesis A}, for $r\in (0,\frac{1}{2})$, we get
\begin{eqnarray*}
\|I_2\|^2_{W^{r,2}([0,T];H)}
&\leq& C(r, T)\sup_{t\in[0,T]}(1+|Y^0|^2)\int^T_0|h(s)|^2_Uds.
\end{eqnarray*}
Thus, we have
\begin{eqnarray*}
\mathbb{E}\|I_2\|_{W^{r,2}([0,T];H)}\leq C_2(r,M,K,T).
\end{eqnarray*}


Collecting the previous inequalities, 
we obtain that for $r\in(0,\frac{1}{2})$,
\begin{eqnarray*}
\mathbb{E}\|R^h\|_{W^{r,2}([0,T];V')}\leq C_{3}(r,M,K,T).
\end{eqnarray*}
In view of Theorem \ref{thm-4}, $R^h$ are uniformly bounded in the space
\[
 \Lambda:=L^2([0,T];D(A))\cap W^{r,2}([0,T];V').
\]
By Lemma \ref{lem-1}, we know $\Lambda$ is compactly imbedded in $L^2([0,T];V)$.
Thus, we obtain $R^h$ is compact in $L^2([0,T];V)$.
\end{flushleft}
 $\hfill\blacksquare$

\subsubsection{Tightness of $\bar{Z}^\varepsilon$}
For any  $h^{\varepsilon}\in \mathcal{A}_M$, consider
\begin{eqnarray}\label{e-47}\notag
d\bar{Z}^\varepsilon(t)&+&A\bar{Z}^\varepsilon(t)+B(\bar{Z}^\varepsilon(t),Y^0+\sqrt{\varepsilon}\lambda(\varepsilon)\bar{Z}^\varepsilon)dt+B(Y^0,\bar{Z}^\varepsilon(t))dt+G(\bar{Z}^\varepsilon(t))dt\\
&=&\lambda^{-1}(\varepsilon)\psi(t,Y^0+\sqrt{\varepsilon}\lambda(\varepsilon)\bar{Z}^\varepsilon)dW(t)+\psi(t,Y^0+\sqrt{\varepsilon}\lambda(\varepsilon)\bar{Z}^\varepsilon)h^{\varepsilon}(t)dt,
\end{eqnarray}
with $\bar{Z}^\varepsilon(0)=0$, then  $\mathcal{G}^{\varepsilon}(W(\cdot)+\lambda(\varepsilon)\int^{\cdot}_0h^{\varepsilon}(s)ds)=\bar{Z}^\varepsilon$.
In view of properties of $Z^\varepsilon$ and $R^h$, by the same method as Theorem \ref{thm-4}, we can obtain
\begin{lemma}\label{prp-4}
Let $Y^0\in \mathcal{H}$. Under \textbf{Hypothesis A} and \textbf{Hypothesis B}, for any family $\{h^\varepsilon, \varepsilon>0\}\subset \mathcal{A}_M$, there exists $\varepsilon_0>0$ such that
\begin{eqnarray}\label{e-49}
\sup_{0\leq \varepsilon\leq \varepsilon_0}\mathbb{E}(\sup_{0\leq t\leq T}|\bar{Z}^\varepsilon|^2+\int^T_0 \|\bar{Z}^\varepsilon\|^2 dt)\leq C(M,K,T).\\
\label{ee-17}
\sup_{0\leq \varepsilon\leq \varepsilon_0}\mathbb{E}(\sup_{0\leq t\leq T}|\partial_z \bar{Z}^\varepsilon|^2+\int^T_0 \|\partial_z\bar{Z}^\varepsilon\|^2 dt)\leq C(M,K,T).
\end{eqnarray}
\end{lemma}

Let $\mathcal{D}(Z)$ be the distribution of $Z$, we have
\begin{prp}\label{prp-6}
Let $Y^0\in \mathcal{H}$. Under \textbf{Hypothesis A} and \textbf{Hypothesis B}, $\mathcal{D}(\bar{Z}^\varepsilon)_{\varepsilon\in[0,\varepsilon_0]}$ is tight in $L^2([0,T];H).$
\end{prp}
\begin{flushleft}
\textbf{Proof.} \quad
Define
\[
F(\bar{Z}^\varepsilon,Y^0):=A\bar{Z}^\varepsilon+B(\bar{Z}^\varepsilon,Y^0+\sqrt{\varepsilon}\lambda(\varepsilon)\bar{Z}^\varepsilon)
+B(Y^0,\bar{Z}^\varepsilon)+G(\bar{Z}^\varepsilon).
\]
From (\ref{e-47}), we have
\begin{eqnarray}\label{e-50}\notag
\bar{Z}^\varepsilon(t)&=&-\int^t_0F(\bar{Z}^\varepsilon(s),Y^0(s))ds\\ \notag
&+&\lambda^{-1}(\varepsilon)\int^t_0\psi(s,Y^0+\sqrt{\varepsilon}\lambda(\varepsilon)\bar{Z}^\varepsilon)dW(s)
+\int^t_0\psi(s,Y^0+\sqrt{\varepsilon}\lambda(\varepsilon)\bar{Z}^\varepsilon)h^{\varepsilon}(s)ds\\
&:=&J_1(t)+J_2(t)+J_3(t).
\end{eqnarray}
Let $\alpha>1$ be fixed, we deduce from the interpolation inequality that $D(A^{\frac{\alpha}{2}})\subset L^{\infty}(\mathcal{M})$. Then, with the aid of the H\"{o}lder inequality and the Cauchy-Schwarz inequality, we have
\begin{eqnarray*}
\|F(\bar{Z}^\varepsilon,Y^0)\|_{D(A^{-\frac{\alpha}{2}})}&\leq& C\|\bar{Z}^\varepsilon\|+C|\bar{Z}^\varepsilon|\|Y^0\|+C\|\bar{Z}^\varepsilon\||\partial_z Y^0|+C\sqrt{\varepsilon}\lambda(\varepsilon)(|\bar{Z}^\varepsilon|+|\partial_z \bar{Z}^\varepsilon|)\|\bar{Z}^\varepsilon\|\\
&&\ +C|Y^0|\|\bar{Z}^\varepsilon\|+C|\partial_z \bar{Z}^\varepsilon|\|Y^0\|+C|\bar{Z}^\varepsilon|.
\end{eqnarray*}

By the H\"{o}lder inequality, we obtain
\begin{eqnarray*}
&&\|J_1(t)-J_1(s)\|^2_{D(A^{-\frac{\alpha}{2}})}\\
&=&\|\int^t_s F(\bar{Z}^\varepsilon,Y^0)dl \|^2_{D(A^{-\frac{\alpha}{2}})}\\
&\leq& C(t-s)\int^t_s\| F(\bar{Z}^\varepsilon,Y^0)\|^2_{D(A^{-\frac{\alpha}{2}})}dl\\
&\leq& C(t-s)\Big[\int^T_0\|\bar{Z}^\varepsilon(t)\|^2dt+C\sup_{t\in [0,T]}|\bar{Z}^\varepsilon|^2\int^T_0\|Y^0(t)\|^2dt\\
&&\ +C\sup_{t\in [0,T]}|\partial_z Y^0|^2\int^T_0\|\bar{Z}^\varepsilon\|^2dt+C\varepsilon\lambda^2(\varepsilon)\sup_{t\in [0,T]}(|\bar{Z}^\varepsilon|^2+|\partial_z \bar{Z}^\varepsilon|^2)\int^T_0\|\bar{Z}^\varepsilon\|^2dt\\
&&\ +C\sup_{t\in [0,T]}|Y^0(t)|^2\int^T_0\|\bar{Z}^\varepsilon(t)\|^2dt+C\sup_{t\in [0,T]}|\partial_z \bar{Z}^\varepsilon(t)|^2 \int^T_0\|Y^0(t)\|^2dt+C(t-s)\sup_{t\in [0,T]}|\bar{Z}^\varepsilon(t)|\Big].
\end{eqnarray*}
From the definition of $W^{r,2}([0,T]; D(A^{-\frac{\alpha}{2}}))$, for $r\in (0,\frac{1}{2})$, we have
\begin{eqnarray*}
&&\|J_1\|^2_{W^{r,2}([0,T];D(A^{-\frac{\alpha}{2}}))}\\
&=&\int^T_0\|J_1(t)\|^2_{D(A^{-\frac{\alpha}{2}})}dt+\int^T_0\int^T_0\frac{\|J_1(t)-J_1(s)\|^2_{D(A^{-\frac{\alpha}{2}})}}{|t-s|^{1+2r}}dsdt\\
&\leq&
C(r, T)\Big[\int^T_0\|\bar{Z}^\varepsilon(t)\|^2dt+C\sup_{t\in [0,T]}|\bar{Z}^\varepsilon|^2\int^T_0\|Y^0(t)\|^2dt\\
&&\ +C\sup_{t\in [0,T]}|\partial_z Y^0|^2\int^T_0\|\bar{Z}^\varepsilon\|^2dt+C\varepsilon\lambda^2(\varepsilon)\sup_{t\in [0,T]}(|\bar{Z}^\varepsilon|^2+|\partial_z \bar{Z}^\varepsilon|^2)\int^T_0\|\bar{Z}^\varepsilon\|^2dt\\
&&\ +C\sup_{t\in [0,T]}|Y^0(t)|^2\int^T_0\|\bar{Z}^\varepsilon(t)\|^2dt+C\sup_{t\in [0,T]}|\partial_z \bar{Z}^\varepsilon(t)|^2 \int^T_0\|Y^0(t)\|^2dt+C(t-s)\sup_{t\in [0,T]}|\bar{Z}^\varepsilon(t)|\Big].
\end{eqnarray*}
Thus, by (\ref{ee-9}) and (\ref{ee-10}), we have
\begin{eqnarray*}
\mathbb{E}\|J_1\|_{W^{r,2}([0,T];D(A^{-\frac{\alpha}{2}}))}
\leq C_{1}(r,M,K,T).
\end{eqnarray*}

Referring to (\ref{ee-9}), (\ref{e-49}) and by \textbf{Hypothesis A}, for any $r\in (0,\frac{1}{2})$ and $p\geq 2$, we have
\begin{eqnarray*}
\mathbb{E}\|J_2\|_{W^{r,2}([0,T]; H)}
&\leq& C(r)\mathbb{E}\int^T_0\|\psi(s,Y^0+\sqrt{\varepsilon}\lambda(\varepsilon)\bar{Z}^\varepsilon)\|^2_{\mathcal{L}_2(U;H)}ds\\
&\leq&C(r,K,T)\mathbb{E}\sup_{t\in[0,T]}(1+|Y^0|^2+\varepsilon\lambda^2(\varepsilon)|\bar{Z}^\varepsilon|^2)\\
&\leq&C_2(r,K,T).
\end{eqnarray*}
By \textbf{Hypothesis A}, for $r\in (0,\frac{1}{2})$, we get
\begin{eqnarray*}
\|J_3\|^2_{W^{r,2}([0,T];H)}
\leq C(r, T)\sup_{t\in[0,T]}(1+|Y^0|^2+\varepsilon\lambda^2(\varepsilon)|\bar{Z}^\varepsilon|^2)\int^T_0|h^\varepsilon(s)|^2_Uds.
\end{eqnarray*}
Thus, from (\ref{ee-9}) and (\ref{e-49}), we have
\begin{eqnarray*}
\mathbb{E}\|J_3\|_{W^{r,2}([0,T];H)}\leq C_3(r,M,K,T).
\end{eqnarray*}


Collecting the previous inequalities, 
we obtain that for $r\in(0,\frac{1}{2})$,
\begin{eqnarray}\label{e-51}
\mathbb{E}\|\bar{Z}^\varepsilon\|_{W^{r,2}([0,T];D(A^{-\frac{\alpha}{2}}))}\leq C_{4}(r,M,K,T).
\end{eqnarray}
In view of Lemma  \ref{prp-4}, $\bar{Z}^\varepsilon$ are bounded uniformly in $\varepsilon$ in the space
\[
 \Lambda:=L^2([0,T];V)\cap W^{r,2}([0,T];D(A^{-\frac{\alpha}{2}})).
\]
By Lemma \ref{lem-1}, we know $\Lambda$ is compactly imbedded in $L^2([0,T];H)$.
Denote
\begin{equation*}
  |\cdot|_{\Lambda}:=|\cdot|_{L^2([0,T];V)}+\|\cdot\|_{W^{r,2}([0,T];D(A^{-\frac{\alpha}{2}}))},
\end{equation*}
then for any $L>0$,
 \[
 K_L=\{Y\in L^2([0,T];H), |Y|_{\Lambda}\leq L\}
 \]
 is compact in $L^2([0,T];H)$.
Finally, notice that
 \[
 P(\bar{Z}^\varepsilon\in K^c_L)\leq P(|\bar{Z}^\varepsilon|_{\Lambda}> L)\leq\frac{E(|\bar{Z}^\varepsilon|_{\Lambda})}{L}\leq \frac{C}{L},
 \]
choosing sufficiently large constant $L$, we obtain $\mathcal{D}(\bar{Z}^\varepsilon)_{\varepsilon\in[0,\varepsilon_0]}$ is tight in $L^2([0,T];H)$.
\end{flushleft}

 $\hfill\blacksquare$

Arguing similarly on the term $J_2$ and by (\ref{ee-9}), we apply Theorem 2.2 in \cite{FG95}, we have the imbedding $C([0,T]; H)\cap W^{r,p}([0,T];D(A^{-\frac{\alpha}{2}}))\subset C([0,T];D(A^{-\frac{\alpha}{2}}))$ is compact. Hence by (\ref{e-51}), we obtain
\begin{prp}\label{prp-7}
$\mathcal{D}(\bar{Z}^\varepsilon)_{\varepsilon\in[0,\varepsilon_0]}$ is tight in $C([0,T];D(A^{-\frac{\alpha}{2}}))$.
\end{prp}

\subsection{Proof of Theorem \ref{thm-8}}
According to Theorem \ref{thm-2}, the proof of Theorem \ref{thm-8} will be completed if the following Theorem \ref{thm-9} and Theorem \ref{thm-1} are established.

\begin{thm}\label{thm-9}
For $Y_0\in V$. Under \textbf{Hypothesis A} and \textbf{Hypothesis B}, for every positive constant $M>0$, the family
\[
K_M=\left\{\mathcal{G}^0(\int^{\cdot}_{0}h(s)ds): h\in T_M\right\}
\]
is a compact subset of $C([0,T];H)\cap L^2([0,T];V)$.
\end{thm}
\begin{flushleft}
\textbf{Proof}.\quad
Let $\{R^h=\mathcal{G}^0(\int^{\cdot}_0 h^n(s)ds); n\geq 1\}$ be a sequence of elements in $K_M$. With the aid of estimates in Theorem \ref{thm-4}, we assert that there exists a subsequence still denoted by $\{n\}$ and $h\in T_M$ such that
\begin{description}
  \item[(i)] $h^{n}\rightarrow h$ in $T_M$ weakly as $n \rightarrow \infty$.
  \item[(ii)] $R^{h^{n}}\rightarrow R^h$ in $L^2([0,T]; D(A))$ weakly.
  \item[(iii)] $R^{h^{n}}\rightarrow R^h$ in $ L^{\infty}([0,T]; V)$ weak-star.
  \item[(iv)] $R^{h^{n}}\rightarrow R^h$ in $L^2([0,T]; V)$ strongly.
\end{description}

Using the similar method as Sect. 5.3 in \cite{R-R}, we have $R^h=\mathcal{G}^0(\int^{\cdot}_0 h(s)ds)$ and $R^h\in C([0,T]; V) \cap L^2([0,T]; D(A))$.

Here, we only prove that $R^{h^{n}}\rightarrow R^h$ in $C([0,T];H)\cap L^2([0,T];V)$. Denote  $\rho^n(t)=R^{h^{n}}-R^h$, using (\ref{equ-9}), we obtain
\begin{eqnarray}\label{e-81}
d\rho^n(t)+A\rho^n(t)dt+B(\rho^n(t), Y^0)dt+B( Y^0,\rho^n(t))dt+G(\rho^n(t))dt=\psi(t, Y^0)(h^n-h)dt.
\end{eqnarray}
Applying the chain rule to (\ref{e-81}), it gives
\begin{eqnarray*}
|\rho^n(t)|^2+2\int^t_0\|\rho^n(s)\|^2ds
&=&-2\int^t_0\langle B(\rho^n(s), Y^0(s)), \rho^n(s) \rangle ds-2\int^t_0\langle B(Y^0(s), \rho^n(s)), \rho^n(s) \rangle ds\\
&&\ -2\int^t_0( G(\rho^n(s)), \rho^n(s)) ds+2\int^t_0(\psi(s, Y^0(s))(h^n-h),\rho^n(s) )ds\\
&\leq& 2\int^t_0|\langle B(\rho^n(s), Y^0(s)), \rho^n(s) \rangle| ds+2\int^t_0|( G(\rho^n(s)), \rho^n(s))| ds\\
&&\ +2\int^t_0|(\psi(s, Y^0(s))(h^n-h),\rho^n(s))|ds.
\end{eqnarray*}
From (\ref{ee-3}) and the H\"{o}lder inequality, we get
\begin{eqnarray*}
&&\sup_{t\in [0,T]}|\rho^n(t)|^2+2\int^T_0\|\rho^n(t)\|^2dt\\
&\leq& C\int^T_0(\|Y^0\||\rho^n|\|\rho^n\|+|\partial_z Y^0|\|\rho^n\||\rho^n|^{\frac{1}{2}}\|\rho^n\|^{\frac{1}{2}}) ds
+C\int^T_0|\rho^n(t)|\|\rho^n(t)\| dt\\
&&\ +C\int^T_0\|\psi(s, Y^0(s)\|_{\mathcal{L}_2(U;H)}|h^n-h|_{U}|\rho^n(s))|ds\\
&\leq& \int^T_0\|\rho^n(t)\|^2dt+C\int^T_0(\|Y^0\|^2+|\partial_z Y^0|^4+1)|\rho^n(s)|^2 ds\\
&&\ +CM^{\frac{1}{2}}(1+\sup_{t\in[0,T]}|Y^0(t)|)\Big(\int^T_0|\rho^n(s)|^2 ds\Big)^{\frac{1}{2}},
\end{eqnarray*}
where
\[
\int^T_0|h^n(t)-h(t)|^2_{U}dt\leq 2\int^T_0|h^n(t)|^2dt+2\int^T_0|h(t)|^2dt\leq 4M
\]
is used. Then
\begin{eqnarray}\notag
&&\sup_{t\in [0,T]}|\rho^n(t)|^2+\int^T_0\|\rho^n(t)\|^2dt\\
\label{e-82}
&\leq& C\int^T_0(\|Y^0\|^2+|\partial_z Y^0|^4+1)|\rho^n(s)|^2 ds +CM^{\frac{1}{2}}(1+\sup_{t\in[0,T]}|Y^0(t)|)\left(\int^T_0|\rho^n(s)|^2 ds\right)^{\frac{1}{2}}.
\end{eqnarray}
Applying Gronwall inequality to (\ref{e-82}), we get
\begin{eqnarray}\notag
&&\sup_{t\in [0,T]}|\rho^n(t)|^2+\int^T_0\|\rho^n(t)\|^2dt\\
\label{e-83}
&\leq& CM^{\frac{1}{2}}(1+\sup_{t\in[0,T]}|Y^0(t)|)\Big(\int^T_0|\rho^n(s)|^2 ds\Big)^{\frac{1}{2}}\times \exp\Big\{\int^T_0(\|Y^0(t)\|^2+|\partial_z Y^0(t)|^4+1)dt\Big\}.
\end{eqnarray}
From \textbf{(iv)}, (\ref{ee-9}) and  (\ref{ee-10}), we have
\begin{eqnarray}
\sup_{t\in [0,T]}|\rho^n(t)|^2+\int^T_0\|\rho^n(t)\|^2dt\rightarrow 0.
\end{eqnarray}
We complete the proof.

\end{flushleft}

$\hfill\blacksquare$

\begin{thm}\label{thm-1}
Fix $M<\infty$ and $h^\varepsilon\subset \mathcal{A}_M$ converges in distribution to $h$ as $\varepsilon\rightarrow 0$. Then
\[
\mathcal{G}^{\varepsilon}\Big(W(\cdot)+\lambda(\varepsilon)\int^{\cdot}_0h^\varepsilon(s)ds\Big)\ converges\ in\ distribution\ to\ \mathcal{G}^{0}\Big(\int^{\cdot}_0h(s)ds\Big),
\]
in $C([0,T];H)\cap L^2([0,T]; V)$ with a good rate function given by (\ref{eq-5}) as $\varepsilon\rightarrow 0$.
\end{thm}
\begin{flushleft}
\textbf{Proof}.\quad Suppose that $\{h^\varepsilon\}_{\varepsilon>0}\subset \mathcal{A}_M$ and $h^\varepsilon$ converges to $h$ as $T_M-$valued random elements in distribution.
By Girsanov's Theorem, we obtain $\bar{Z}^\varepsilon=\mathcal{G}^{\varepsilon}(W(\cdot)+\lambda(\varepsilon)\int^{\cdot}_0h^\varepsilon(s)ds)$.
Consider
\[
dX^\varepsilon(t)+AX^\varepsilon(t)dt=\lambda^{-1}(\varepsilon)\psi(t, Y^0(t)+\sqrt{\varepsilon}\lambda(\varepsilon)X^\varepsilon(t))dW(t),
\]
with initial value $X^\varepsilon(0)=0$.
Applying It\^{o} formula, \textbf{Hypothesis A}, (\ref{ee-9}), we get
\begin{eqnarray}\label{e-84}
\lim_{\varepsilon\rightarrow 0}\Big[\mathbb{E}\sup_{t\in [0,T]}\|X^\varepsilon(t)\|^2+\mathbb{E}\int^T_0\|X^\varepsilon(t)\|^2_{D(A)}dt\Big]=0.
\end{eqnarray}
Using the same method as Proposition \ref{prp-6-1}, we deduce that $X^\varepsilon$ is tight in $C([0,T];H)\cap L^2(0,T;V)$.

Set
\[
\prod=\left(L^2([0,T];H)\cap C([0,T];D(A^{-\frac{\alpha}{2}})), T_M, C([0,T];H)\cap L^2([0,T];V)\right).
\]
By Proposition \ref{prp-6} and Proposition \ref{prp-7}, we know that the family $\{(\bar{Z}^\varepsilon, h^\varepsilon, X^\varepsilon), \varepsilon \in (0,\varepsilon_0)\}$ is tight in $\prod$.  Let $(Z, h, 0)$ be any limit point of $\{(\bar{Z}^\varepsilon, h^\varepsilon, X^\varepsilon), \varepsilon \in (0,\varepsilon_0)\}$. We will show that $Z$ has the same law as $\mathcal{G}^0(\int^{\cdot}_0 h(s)ds)$, and in fact $\bar{Z}^\varepsilon$ converges in distribution to $Z$ in $C([0,T];H)\cap L^2([0,T]; V)$ as $\varepsilon\rightarrow 0$, which implies Theorem \ref{thm-1}.

By the Skorokhod representative theorem, there exists a stochastic basis $(\Omega^1, \mathcal{F}^1, \{\mathcal{F}^1_t\}_{t\in [0,T]}, {P}^1)$ and, on this basis, $\prod-$valued random variables $(\tilde{U}^{\varepsilon}, \tilde{h}^{\varepsilon},\tilde{X}^{\varepsilon}),(\tilde{U}, \tilde{h},0)$  such that $(\tilde{U}^{\varepsilon}, \tilde{h}^{\varepsilon},\tilde{X}^{\varepsilon})$ (respectively $(\tilde{U}, \tilde{h},0)$) has the same law as $(\bar{Z}^\varepsilon, h^\varepsilon, X^\varepsilon)$ (respectively $(Z,h,0)$), and $(\tilde{U}^{\varepsilon}, \tilde{h}^{\varepsilon},\tilde{X}^{\varepsilon})\rightarrow (\tilde{U}, \tilde{h},0)$, ${P}^1-$a.s. in $\prod$. From the equation satisfied by $(\bar{Z}^\varepsilon, h^\varepsilon, X^\varepsilon)$, we see that $(\tilde{U}^{\varepsilon}, \tilde{h}^{\varepsilon},\tilde{X}^{\varepsilon})$ satisfies the following equation in the distribution sense as follows:
\begin{eqnarray}\notag
d(\tilde{U}^{\varepsilon}(t)-\tilde{X}^{\varepsilon})&+&A(\tilde{U}^{\varepsilon}(t)-\tilde{X}^{\varepsilon})dt +(B(\tilde{U}^{\varepsilon}(t),Y^0(t)+\sqrt{\varepsilon}\lambda(\varepsilon)\tilde{X}^{\varepsilon}(t))+B(Y^0(t),\tilde{U}^{\varepsilon}(t)))dt\\
\label{e-60}
&=&\psi(t,Y^0(t)+\sqrt{\varepsilon}\lambda(\varepsilon)\tilde{U}^{\varepsilon})\tilde{h}^{\varepsilon}(t),
\end{eqnarray}
and
\begin{eqnarray*}
&&P^1(\tilde{U}^{\varepsilon}-\tilde{X}^{\varepsilon}\in C([0,T];H)\cap L^2([0,T];V))\\
&=& P(\bar{Z}^{\varepsilon}-{X}^{\varepsilon}\in C([0,T];H)\cap L^2([0,T];V))\\
&=&1
\end{eqnarray*}
Let ${\Omega}^1_0$ be the subset of $\Omega^1$ such that $(\tilde{U}^{\varepsilon}, \tilde{h}^{\varepsilon},\tilde{X}^{\varepsilon})\rightarrow (\tilde{U}, \tilde{h},0)$ in $\prod$, we have $P^1({\Omega}^1_0)=1$. For any $\tilde{\omega}\in {\Omega}^1_0$, we have
\begin{eqnarray}\label{e-89}
\sup_{t\in [0,T]}|\tilde{U}^{\varepsilon}(\tilde{\omega},t)-\tilde{U}(\tilde{\omega},t)|^2
+\int^T_0\|\tilde{U}^{\varepsilon}(\tilde{\omega},t)-\tilde{U}(\tilde{\omega},t)\|^2dt\rightarrow 0 \quad as\ \varepsilon \rightarrow 0.
\end{eqnarray}
Set $\tilde{\eta}^\varepsilon=\tilde{U}^{\varepsilon}-\tilde{X}^{\varepsilon}$, then $\tilde{\eta}^\varepsilon(\tilde{\omega})\in C([0,T];H)\cap L^2([0,T];V)$, and $\tilde{\eta}^\varepsilon(\tilde{\omega})$ satisfies
\begin{eqnarray}\notag
d\tilde{\eta}^\varepsilon(t)&+&A\tilde{\eta}^\varepsilon(t)dt
+\Big[B\Big(\tilde{\eta}^\varepsilon(t)+\tilde{X}^{\varepsilon}(t), Y^0(t)+\sqrt{\varepsilon}\lambda(\varepsilon)(\tilde{\eta}^\varepsilon(t)+\tilde{X}^{\varepsilon}(t))\Big)
+B\Big(Y^0(t),\tilde{\eta}^\varepsilon(t)+\tilde{X}^{\varepsilon}(t)\Big)\Big]dt\\
\label{e-87}
&=&\psi\Big(t,Y^0(t)+\sqrt{\varepsilon}\lambda(\varepsilon)(\tilde{\eta}^\varepsilon(t)+\tilde{X}^{\varepsilon}(t))\Big)\tilde{h}^{\varepsilon}(t)dt,
\end{eqnarray}
with initial value $\tilde{\eta}^\varepsilon(0)=0$.

Consider $\hat{U}(t)$ satisfies
\begin{eqnarray}\label{e-88}
\hat{U}(t)=-\int^t_0A\hat{U}(s)ds-\int^t_0(B(\hat{U}(s), Y^0(s))+B(Y^0(s),\hat{U}(s)))ds+\int^t_0\psi(s, Y^0(s))\tilde{h}(s)ds.
\end{eqnarray}
Since
\begin{eqnarray*}
\lim_{\varepsilon\rightarrow 0}\Big[\sup_{t\in [0,T]}|\tilde{X}^\varepsilon(\tilde{\omega},t)|^2+\int^T_0\|\tilde{X}^\varepsilon(\tilde{\omega},t)\|^2dt\Big]=0, \quad\tilde{U}^{\varepsilon}=\tilde{\eta}^\varepsilon+\tilde{X}^{\varepsilon},
\end{eqnarray*}
using (\ref{e-87}), (\ref{e-88}) and by the same method as Theorem \ref{thm-9}, we obtain
\begin{eqnarray}\label{e-86}
\lim_{\varepsilon\rightarrow 0}\Big[\sup_{t\in [0,T]}|\tilde{U}^\varepsilon(\tilde{\omega},t)-\hat{U}(\tilde{\omega},t)|^2+\int^T_0\|\tilde{U}^\varepsilon(\tilde{\omega},t)-\hat{U}(\tilde{\omega},t)\|^2dt\Big]=0.
\end{eqnarray}
Hence, from  (\ref{e-89}) and (\ref{e-86}), we deduce that $\tilde{U}=\hat{U}=\mathcal{G}^0(\int^{\cdot}_0 \tilde{h}(s)ds)$, and $\tilde{U}$ has the same law as $\mathcal{G}^0(\int^{\cdot}_0 {h}(s)ds)$. Taking into account that $\bar{Z}^\varepsilon$ and $\tilde{U}^\varepsilon$ has the same law on $C([0,T];H)\cap L^2([0,T];V)$ and by (\ref{e-86}), we deduce that $\mathcal{G}^{\varepsilon}\Big(W(\cdot)+\lambda(\varepsilon)\int^{\cdot}_0h^\varepsilon(s)ds\Big)\ $converges in distribution to $\mathcal{G}^{0}\Big(\int^{\cdot}_0h(s)ds\Big)$ as $\varepsilon \rightarrow 0$. We complete the proof.
\end{flushleft}

$\hfill\blacksquare$

%





\def\refname{ References}

\end{document}